\documentclass[12pt]{article}
\usepackage{amsfonts,amsmath,amssymb}
\usepackage{enumerate}
\usepackage{theorem}
\numberwithin{equation}{section}

\newtheorem{theorem}{Theorem}[section]
\newtheorem{lemma}[theorem]{Lemma}
\newtheorem{proposition}[theorem]{Proposition}

\newtheorem{corollary}[theorem]{Corollary}
\newtheorem{remark}[theorem]{Remark}

\newcommand{\cl}[1]{\mathcal{#1}} 

\begin{document}

\title{Decompositions of Reflexive Bimodules over Maximal Abelian Selfadjoint Algebras}
\author{G. Eleftherakis}
\date{}
\maketitle

\begin{abstract}
We generalize the notion of `diagonal' from the class of CSL
algebras to masa bimodules. We prove that a reflexive masa
bimodule decomposes as a sum of two bimodules, the diagonal and a
module generalizing the w*-closure of the Jacobson radical of a
CSL algebra. The latter module turns out to be reflexive, a result
which is new even for CSL algebras. We show that the projection
onto the direct summand contained in the diagonal is contractive
and preserves compactness and reduces rank of operators. Stronger
results are obtained when the module is the reflexive hull of its
rank-one subspace.

\medskip



\end{abstract}

\section{Introduction}
In this paper we attempt a generalisation of the concept of the
diagonal of a CSL algebra to reflexive spaces of operators which
are modules over maximal abelian selfadjoint algebras (masas).

Recall \cite{dav} that a CSL algebra is an algebra $\cl A$ of
operators on a Hilbert space $H$ which can be written in the form
 $$\cl A = \{A\in B(H) : AP=PAP \; \text{for all }P\in \cl S\}$$
 where $\cl S$ is a commuting family of projections. Note that
 $\cl A$ contains any masa containing $\cl S''$.

More generally, a reflexive masa bimodule $\cl U$ of operators
from $H$ to another Hilbert space $K$ can be written in the form
 $$\cl U = \{T\in B(H,K) : TP=\phi(P)TP \; \text{for all }P\in \cl S\}$$
 where $\cl S$ is a commuting family of projections on $H$ and
 $\phi$ maps them to commuting projections on $K$ (see below for
 details).

 The diagonal $\mathcal{A\cap A^*}$ of a CSL
algebra $\cl A$ is a von Neumann algebra, which equals the
commutant
 $$\cl S'=\{A\in B(H) : AP=PA \; \text{for all }P\in \cl S\}$$
of the corresponding invariant projection family. The natural
corresponding object for a reflexive masa bimodule $\cl{U}$ is a
ternary ring of operators (TRO)
$$\Delta (\cl{U})=\{T\in B(H,K) : TP=\phi(P)T \; \text{for all }P\in \cl S\}$$
which is also a reflexive masa bimodule.

This `diagonal' $\Delta (\cl{U})$ is the primary object of study
of the present paper.

We decompose $\cl{U}$ as a sum $\cl{U}_0+\Delta (\cl{U})$, where
$\cl{U}_0$ also turns out to be reflexive (Theorem 5.2). This is
new even for the case of CSL algebras; note, however, that for
nest algebras reflexivity of w*-closed bimodules is automatic
\cite{erpo}. An analogous decomposition for the case of nest
subalgebras of von Neumann algebras is in \cite{ls}.\\ We also
prove (Corollary 5.3) that the bimodule $\cl{U}_0$ has in our
context the role corresponding to the $w^*$ closure of the
Jacobson radical of a CSL algebra.

The diagonal $\Delta (\cl{U})$ is proved to be generated by a partial isometry and natural von Neumann algebras assosiated to $\cl{U}$ (Theorem 4.1).

The above decomposition may be further refined to a direct sum:
$\cl U=\cl{U}_0\oplus \cl{M}$ where $\cl{M}$ is a TRO ideal of the
diagonal $\Delta (\cl{U})$ (Theorem 3.3), containing
the compact operators of the diagonal (Proposition 6.3).

In case $\cl{U}$ is strongly reflexive (that is, coincides with
the reflexive hull of the rank one operators it contains) we show (Theorem 7.4)
that $\cl{M}$ coincides with the w*-closed linear span of the
finite rank operators of the diagonal, an equality which fails in
general.

As in the case of von Neumann algebras, we show that every TRO
decomposes in an `atomic' and a \lq nonatomic' part. The `atomic'
part of the diagonal $\Delta (\cl{U})$ is contained (properly in
general) in $\cl{M}$ (Proposition 6.3).

We also study the projection $\theta :\cl{U}\longrightarrow \cl{M}$
defined by the above direct sum decomposition. We prove that it is
contractive and maps compact operators to compact operators and 
finite rank operators to operators of at most the same rank.

In case $\cl{U}$ is strongly reflexive, we show that $\theta
=D|_{\cl{U}}$, where $D$ is the natural projection onto the `atomic' part of the diagonal $\Delta (\cl U)$.

A main tool used to obtain these results is an appropriate
sequence of projections $(U_n)$ on $B(H,K)$  which depend on
$\cl{U}$ . This sequence behaves analogously to the net of
`diagonal sums' used in nest algebras (see for example
\cite{dav}).

In nest algebra theory, the net of diagonal sums of a compact
operator converges in norm to a compact operator in the `atomic'
part of the diagonal. This has been generalised to CSL algebras by
Katsoulis \cite{kats}. Here we show (Proposition 6.10) that for every compact
operator $K$, the sequence $(U_n(K))$ converges in norm to $D(K)$.

\bigskip

 We present some definitions and concepts we use in this work. All Hilbert spaces will be assumed separable.
\medskip

If $S$ is a set of operators then $R_1(S)$ denotes the subset of $S$
which contains the rank 1 operators and the zero operator. If $H$
is a Hilbert space and $S\subset B(H),$ the set of orthogonal
projections of $S$ is denoted by $\cl{P}(S).$

\medskip

If $H_1,H_2$ are Hilbert spaces, $C_1(H_1,H_2)$ are the trace class operators and $\cl{R}$ a subset of $C_1(H_1,H_2),$ we denote by $\cl{R}^0$ the set of operators which are annihilated by $\cl{R}:$
$$\cl{R}^0=\{T\in B(H_2,H_1): tr(TS)=0\;\text{for\;all}\;S\in \cl{R}\}.$$

\bigskip

Let $H_1,H_2$ be Hilbert spaces and $\cl{U}$\;a subset of
$B(H_1,H_2).$ Then the \textbf{reflexive hull} of $\cl{U}$ is
defined \cite{logshu} to be the space
 $$\mathrm{Ref} (\cl{U})=\{T\in
B(H_1,H_2):Tx\in\overline{[\cl{U}x]}\;\text{for
each}\;x\in\;H_1\}.$$
Simple arguments show that  $$\mathrm{Ref}
(\cl{U})=\{T\in B(H_1,H_2):\text{for all
projections}\;E,F:E\cl{U}F=0\Rightarrow ETF=0\}$$A subspace
$\cl{U}$ is called \textbf{reflexive} if  $\cl{U}=\mathrm{Ref} (\cl{U}).$
It is called \textbf{strongly reflexive} if there exists a set
$L\subset B(H_1,H_2)$\;of rank $1$\;operators such
that\;$\cl{U}=\mathrm{Ref} (L).$

\bigskip

Now we present some concepts introduced by Erdos \cite{erd}.

\medskip

Let $\cl{P}_i=\cl{P}(B(H_i)), i=1,2.$ Define $\phi=\mathrm{Map} (\cl{U})$ to be the map $\phi:\cl{P}_1\rightarrow \cl{P}_2$\;which associates to every
$P\in\cl{P}_1$\;the projection onto the subspace
$[TPy:T\in\;\cl{U},y\in\;H_1]^{-}.$ The map $\phi$\;is $\vee-$continuous (that is, it preserves arbitrary suprema) and $0$ preserving.

\medskip

Let
$\phi^*=\mathrm{Map} (\mathcal{U}^*),\mathcal{S}_{1,\phi}=\{\phi^*(P)^{\bot}:P\in\cl{P}_2\},\mathcal{S}_{2,\phi}=\{\phi(P):P\in\cl{P}_1\}.$
Erdos has proved that $\mathcal{S}_{1,\phi}$\;is meet complete and
contains the identity projection, $\mathcal{S}_{2,\phi}$\;is join
complete and contains the zero projection, while
$\phi|_{\mathcal{S}_{1,\phi}}:\mathcal{S}_{1,\phi}\rightarrow\mathcal{S}_{2,\phi}$\;is
a bijection. In fact
\begin{equation}\label{sx}
(\phi |_{\cl{S}_{1,\phi }})^{-1}(Q)=\phi ^*(Q^\bot )^\bot
\end{equation}
for all $Q\in \cl{S}_{2,\phi }$ and
$$\mathrm{Ref} (\cl{U})=\{T\in\;B(H_1,H_2):\phi(P)^{\bot}TP=0\; \text{for
each}\;P\in\mathcal{S}_{1,\phi}\}.$$
We call the families
$\cl{S}_{1,\phi }, \cl{S}_{2,\phi }$ the \textbf{semilattices} of
$\cl{U}.$

\medskip

A \textbf{C.S.L.} is a complete abelian lattice of projections which contains the identity and the zero projection.

\bigskip

If $\cl{A}_1\subset B(H_1),$ and $\cl{A}_2\subset B(H_2)$\;are algebras, a subspace $\cl{U}\subset B(H_1,H_2)$\;is called an $\cl{A}_1$,$\cl{A}_2-$\textbf{bimodule} if $\cl{A}_2\cl{U}\cl{A}_1\subset\cl{U}.$

\bigskip

A subspace $\cl{M}$\;of $B(H_1,H_2)$\;is called a \textbf{ternary ring of operators} (TRO)\;if $\cl{M}\cl{M}^*\cl{M}\subset\cl{M}.$ Katavolos and Todorov \cite{kt} have proved that a TRO $\cl{M}$ is $w^*$ closed if and only if it is $wot$ closed if and only if it is reflexive. In this case, if $\chi=\mathrm{Map} (\cl{M}),$\;then $$\cl{M}=\{T\in\;B(H_1,H_2):TP=\chi(P)T\;\text{for all}\; P\in\mathcal{S}_{1,\chi}\}.$$

\medskip

They also proved that if $\cl{M}$ is a strongly reflexive TRO, then there exist  families of mutually orthogonal projections $(F_n), (E_n)$ such that $\cl{M}=\sum_{n=1}^\infty \oplus E_nB(H_1,H_2)F_n.$ We present a new proof of this result in Corollary 6.9.

\medskip

The following proposition is easily proved.

\begin{proposition}
Let $H_1,H_2$ be Hilbert spaces, $\cl{A}_1\subset\;B(H_1),\cl{A}_2\subset\;B(H_2)$ masas and $\cl{U}$ a $\cl{A}_1$,$\cl{A}_2-$bimodule. Then
$$\mathrm{Ref}(\cl{U}) =\{T\in\;B(H_1,H_2): E\in\cl{P}(\cl{A}_2),\:F\in\cl{P}(\cl{A}_1),\:E\cl{U}F=0\Rightarrow ETF=0\}.$$
\end{proposition}

The next section contains some preliminary results.

\section{Decomposition of a reflexive TRO.}

In this section we show that a $w^*-$closed TRO decomposes into a
\lq nonatomic\rq\; and a \lq totaly atomic\rq\; part.

\medskip

Let $H_1,H_2$ be Hilbert spaces,\;$\mathcal{M}\subset B(H_1,H_2)$
be a $w^*$-closed TRO and $\cl{B}_1=(\cl{M}^*\cl{M})^{\prime
\prime},\cl{B}_2=(\cl{M}\cl{M}^*)^{\prime \prime}.$

\begin{remark}
We suppose that $\cl{M}_0$ is a $w^*-$closed TRO ideal of
$\cl{M};$ namely, $\cl{M}_0$ is a linear subspace of $\cl{M}$ and
$$\cl{M}_0\cl{M}^*\cl{M}\subset \cl{M}_0,\; \cl{M}\cl{M}^*\cl{M}_0\subset \cl{M}_0.$$ It follows that $\cl{M}\cl{M}_0^*\cl{M}\subset \cl{M}_0$ \cite{eor}.\\
Now, we observe that there exist
projections $Q_i$ in the centre of $\cl{B}_i,\;i=1,2$ such that
$\cl{M}_0=\cl{M}Q_1=Q_2\cl{M}.$ Hence $\cl{M}_0$ is a $\cl{B}_1,
\cl{B}_2-$bimodule.
\end{remark}

\textbf{Proof}

Let $\cl{J}_1=[\cl{M}_0^*\cl{M}_0]^{-w^*}$ and
$\cl{J}_2=[\cl{M}_0\cl{M}_0^*]^{-w^*}.$

\medskip

We can easily verify that $\cl{J}_i$ is an ideal of $\cl{B}_i,
i=1,2.$ Hence there is a projection $Q_i$ in the centre of
$\cl{B}_i$ so that $\cl{J}_i=\cl{B}_iQ_i, i=1,2.$

One easily checks that $$\cl{M}\cl{B}_1\subset \cl{M},\;\;
\cl{B}_2\cl{M}\subset \cl{M},$$
$$\cl{M}\cl{J}_1\subset \cl{M}_0,\;\; \cl{J}_2\cl{M}\subset \cl{M}_0$$

\medskip

We observe that $\cl{M}Q_1\subset \cl{M}\cl{J}_1\subset \cl{M}_0$.

\medskip

Let $T\in \cl{M}_0$ then $T^*T\in \cl{J}_1,$ so $T^*T=T^*TQ_1$ and
thus $T=TQ_1.$ Hence $T\in \cl{M}Q_1.$ We conclude that
$\cl{M}_0\subset \cl{M}Q_1$ and hence equality holds.

\medskip

Similarly one shows that $\cl{M}_0=Q_2\cl{M}. \qquad \Box$

\bigskip

Since $[R_1(\mathcal{M})]^{-w^*}$ is a strongly reflexive TRO, by Proposition 3.5 in \cite{kt} there exist mutually orthogonal projections $(F_n)$ in the centre of $\cl{B}_1$ and $(E_n)$ in the centre of $\cl{B}_2$ such that $[R_1(\cl{M})]^{-w^*}=\sum_{n=1}^\infty \oplus E_nB(H_1,H_2)F_n.$ We write $E=\vee _nE_n, F=\vee _nF_n.$

\begin{theorem}The space $\cl{M}$ decomposes in the following direct sum
$$\mathcal{M}=(\mathcal{M}\cap(R_1(\mathcal{M})^*)^0)\oplus [R_1(\mathcal{M})]^{-w^*}.$$
The spaces $\mathcal{M}\cap(R_1(\mathcal{M})^*)^0$ and $[R_1(\mathcal{M})]^{-w^*}$ are TRO ideals of $\cl{M}.$ Moreover $$[R_1(\cl{M})]^{-w^*}=\cl{M}F=E\cl{M}=E\cl{M}F$$
$$\cl{M}\cap (R_1(\mathcal{M})^*)^0=\cl{M}F^\bot =E^\bot \cl{M}=E^\bot \cl{M}F^\bot.$$

\end{theorem}
\textbf{Proof}

\medskip

We observe that $[R_1(\cl{M})]^{-w^*}$ is a TRO ideal of $\cl{M}.$

\medskip

By Remark 2.1 there exists projection $Q$ in the centre of
$\cl{B}_1$ such that \\$[R_1(\cl{M})]^{-w^*}=\cl{M}Q.$

\medskip

For every $m\in \mathbb{N},$ we have $E_mB(H_1,H_2)F_m\subset \cl{M}Q.$\\It follows that  $E_mB(H_1,H_2)F_m=E_mB(H_1,H_2)F_mQ,$ so $F_m=F_mQ.$ We conclude that $\vee _mF_m=F\leq Q.$\\
Since $F\in \cl{B}_1$ we get $\cl{M}F\subset \cl{M},$ therefore $\cl{M}F=\cl{M}FQ\subset \cl{M}Q.$\\
It follows that $$[R_1(\cl{M})]^{-w^*}=\cl{M}Q\supset \cl{M}F\supset [R_1(\cl{M})]^{-w^*}F=[R_1(\cl{M})]^{-w^*}.$$ We proved that $[R_1(\cl{M})]^{-w^*}=\cl{M}F.$

\medskip

If $M\in \cl{M}$ and $R\in R_1(\cl{M}),$ then $R=RF$ so\; $tr(MF^\bot R^*)\\ =tr(M(RF^\bot )^*)=tr(M0)=0.$

\medskip

We conclude that $$\cl{M}F^\bot \subset \cl{M}\cap (R_1(\mathcal{M})^*)^0$$
Hence $\cl{M}=\cl{M}F^\bot +\cl{M}F\subset \mathcal{M}\cap(R_1(\mathcal{M})^*)^0 + [R_1(\mathcal{M})]^{-w^*}\subset \cl{M}.$

\medskip

It follows that
$$\mathcal{M}=(\mathcal{M}\cap(R_1(\mathcal{M})^*)^0) +
[R_1(\mathcal{M})]^{-w^*}.$$

\medskip

We shall prove that this sum is direct.

\medskip

If $T \in [R_1(\mathcal{M})]^{-w^*}\cap(R_1(\mathcal{M})^*)^0$ then $T=\sum_{n=1}^\infty E_nTF_n.$ If $R$ is a rank $1$ operator then $tr(TR)=\sum_{n=1}^\infty tr(E_nTF_nR)=\sum_{n=1}^\infty tr(TF_nRE_n).$

\medskip

But for every $n\in \mathbb{N},\; tr(TF_nRE_n)=tr(T(E_nR^*F_n)^*)=0$\; since \;$E_nR^*F_n\in R_1(\cl{M})$ and $T\in (R_1(\mathcal{M})^*)^0.$

\medskip

Thus $tr(TR)=0$ for every rank 1 operator R, hence $T=0.$ This shows that $[R_1(\mathcal{M})]^{-w^*}\cap(R_1(\mathcal{M})^\bot)^*=0.$

\medskip

We have shown that  $\mathcal{M}=(\mathcal{M}\cap(R_1(\mathcal{M})^*)^0)\oplus [R_1(\mathcal{M})]^{-w^*}.$

\medskip

Since $\cl{M}=\cl{M}F^\bot \oplus \cl{M}F,\;
[R_1(\cl{M})]^{-w^*}=\cl{M}F$ and \\ $\cl{M}F^\bot \subset
\cl{M}\cap (R_1(\mathcal{M})^*)^0$ we conclude that
$$\cl{M}F^\bot = \cl{M}\cap (R_1(\mathcal{M})^*)^0.$$

\medskip

The equalities $E^\bot \cl{M} = \cl{M}\cap (R_1(\mathcal{M})^*)^0, E\cl{M}=[R_1(\cl{M})]^{-w^*}$ are proved similarly.$\qquad\Box$

\begin{proposition}Let $\theta  :\cl{M}\longrightarrow \cl{M}$ be the projection onto $[R_1(\mathcal{M})]^{-w^*}$ defined by the decomposition in Theorem 2.2. Then $\theta  (T)=\sum_{n=1}^\infty E_nTF_n$ for every $T\in \cl{M}.$
\end{proposition}
\textbf{Proof}

\medskip

Since $\cl{M}$ decomposes as the direct sum of the $\cl{B}_1,
\cl{B}_2-$bimodules \linebreak $\cl{M}\cap (R_1(\cl{M})^*)^0$
and $[R_1(\mathcal{M})]^{-w^*},$\; $\theta  $ is a $\cl{B}_1,
\cl{B}_2-$bimodule map: $$\theta  (B_2TB_1)=B_2\theta  (T)B_1$$ for
every $T\in \cl{M}, B_1\in \cl{B}_1, B_2\in \cl{B}_2.$

\medskip

Since $(E_n)\subset \cl{B}_1, (F_n)\subset \cl{B}_2$ we have that:
$$\theta  (T)=\sum_{n=1}^\infty E_n\theta  (T)F_n=\sum_{n=1}^\infty \theta  (E_nTF_n)=\sum_{n=1}^\infty E_nTF_n.  \qquad\Box$$

\section{Decomposition of a reflexive masa bimodule}

Let $H_1$, $H_2$ be Hilbert spaces,\;$\cl{P}_i=\cl{P}(B(H_i)), i=1,2,$\;$ \mathcal{D}_i\subset
B(H_i),i=1,2$ be masas, $\mathcal{U}\subset B(H_1,H_2)$ be a
reflexive $\mathcal{D}_1$,$\mathcal{D}_2-$bimodule. Write
\begin{align*}
& \phi=\mathrm{Map} (\mathcal{U}),\qquad \phi^*=\mathrm{Map} (\mathcal{U}^*),\\
&\mathcal{S}_{2,\phi}=\phi(\cl{P}_1), \quad  \mathcal{S}_{1,\phi}=\{P^\bot:P\in\phi^*(\cl{P}_2)\}\\
& \cl{A}_2=(\mathcal{S}_{2,\phi})^\prime,\qquad
\cl{A}_1=(\mathcal{S}_{1,\phi})^\prime.
\end{align*}
Observe that $\mathcal{S}_{i,\phi}\subset \;\cl{D}_i$ hence
$\cl{D}_i\subset \;\cl{A}_i,i=1,2.$ We define  $$\cl{U}_0=[\phi
(P)TP^\bot:T\in\;\cl{U},P\in\;\cl{S}_{1,\phi}]^{-w^*},$$ $$\Delta
(\cl{U})=\{T: TP=\phi (P)T \;\text{for all}\;
P\in\;\cl{S}_{1,\phi}\}.$$

\medskip

We remark that $\cl{U}_0$ and $\Delta (\cl{U})$ are $\mathcal{D}_1$,$\mathcal{D}_2-$bimodules contained in $\cl{U}$ and $\Delta(\mathcal{U})$ is a reflexive TRO. We call $\Delta(\mathcal{U})$ the \textbf{diagonal} of $\cl{U}.$

\begin{theorem} $\mathcal{U}=\mathcal{U}_0+\Delta(\mathcal{U}).$
\end{theorem}
\textbf{Proof}

\medskip

As noted in the introduction $$\cl{U}=\{T\in\;B(H_1,H_2):\phi(P)^{\bot}TP=0\; \text{for all}\;P\in\mathcal{S}_{1,\phi}\}.$$

\medskip

Since the Hilbert spaces $H_1,H_2$ are separable we can choose a
sequence $(P_n)\subset \;\cl{S}_{1,\phi}$ such that
$$\cl{U}=\{T\in\;B(H_1,H_2):\phi(P_n)^{\bot}TP_n=0\; \text{for
all}\; n\in\;\mathbb{N} \}.$$

\medskip

We define $$V_n:B(H_1,H_2)\rightarrow  B(H_1,H_2):V_n(T)=\phi (P_n)TP_n+\phi (P_n)^\bot  TP_n^\bot ,n\in \mathbb{N}.$$ One easily checks that $V_n$ is idempotent and a norm contraction.

\medskip

We also define $U_n=V_n\circ V_{n-1}\circ ...\circ V_1, n\in \mathbb{N}.$

\medskip

Let $T\in\;\cl{U},$ then $$T=U_1(T)+\phi(P_1)TP_1^{\bot}$$\;$$U_1(T)=U_2(T)+\phi(P_2)U_1(T)P_2^{\bot}$$ by induction  $$U_{n-1}(T)=U_n(T)+\phi(P_n)U_{n-1}(T)P_n^{\bot}$$ for all $n\in\;\mathbb{N} .$

\medskip

Adding the previous equalities we obtain $$T=U_n(T)+M_n$$ where $$M_n=\phi(P_1)TP_1^{\bot}+\phi(P_2)U_1(T)P_2^{\bot}+...+\phi(P_n)U_{n-1}(T)P_n^{\bot} \in\;\cl{U}_0$$ for all $n\in\;\mathbb{N} .$

\medskip

We observe that $\phi(P_i)^{\bot}U_n(T)P_i=\phi(P_i)U_n(T)P_i^{\bot}=0$ for $i=1,2,...n$ and $\|U_n(T)\|\leq \|U_{n-1}(T)\|\leq ...\leq \|T\|$ for all $n\in\;\mathbb{N} .$

\medskip

The sequence $(U_n(T))$ is bounded, so there exists a subsequence $(U_{n_m}(T))$ that converges in the weak$-*$ topology to an operator L.

\medskip

Then
$M_{n_m}=T-U_{n_m}(T)\stackrel{w^*}{\rightarrow}T-L=M\in \cl{U}_0.$

\medskip

 Since $\phi(P_i)^{\bot}LP_i=\phi(P_i)LP_i^{\bot}=0$
for all $i\in\;\mathbb{N} $ we have $L\in\;\Delta(\mathcal{U})$
and $T=M+L\; \in\; \cl{U}_0+\Delta (\cl{U}).\qquad\Box$

\begin{remark}The following are equivalent:\\
i)\;$\cl{U}$ is a TRO\\ii)\;\;$\cl{U}=\Delta (\cl{U})$\\iii)\;\;\;$\cl{U}_0=0.$
\end{remark}

\begin{theorem}There exist projections $Q_i\in \cl{D}_i, i=1,2$ such that:$$\mathcal{U}=\mathcal{U}_0\oplus (I-Q_2)\Delta
(\mathcal{U})(I-Q_1)=\mathcal{U}_0\oplus (I-Q_2)\Delta
(\mathcal{U})=\mathcal{U}_0\oplus \Delta
(\mathcal{U})(I-Q_1).$$
\end{theorem}
\textbf{Proof}

\medskip

We make the following observations:

\medskip

$i)$ $\cl{U}\Delta (\cl{U})^*\Delta (\cl{U})\subset \cl{U},\;\; \Delta (\cl{U})^*\Delta (\cl{U})\cl{U}\subset \cl{U}.$\\
\emph{Proof}\\
Let $T\in \cl{U}, M,N\in \Delta (\cl{U}).$ Then for every $P\in \cl{S}_{1,\phi }$ we have\\
$\phi (P)^\bot TM^*NP=\phi (P)^\bot TM^*\phi (P)N=\phi (P)^\bot TPM^*N=0M^*N=0.$\\
Thus $TM^*N\in \cl{U}.$ Similarly we have that $MN^*T\in \cl{U}.$

\medskip

$ii)$ $\cl{U}_0\Delta (\cl{U})^*\Delta (\cl{U})\subset \cl{U}_0,\;\; \Delta (\cl{U})\Delta (\cl{U})^*\cl{U}_0\subset \cl{U}_0.$\\
\emph{Proof}\\
Let $T\in \cl{U}, M,N\in \Delta (\cl{U}).$ Then for every $P\in \cl{S}_{1,\phi }$ we have\\
$\phi (P)TP^\bot M^*N=\phi (P)TM^*\phi (P)^\bot N=\phi (P)TM^*NP^\bot .$\\
It follows by $(i)$ that $TM^*N\in \cl{U}$ so $\phi (P)TP^\bot M^*N\in \cl{U}_0.$\\
Taking the $w^*$ closed linear span we get $SM^*N\in \cl{U}_0$ for all $S\in \cl{U}_0, M,N\in \Delta (\cl{U}).$ Similarly we have that $\Delta (\cl{U})\Delta (\cl{U})^*\cl{U}_0\subset \cl{U}_0.$

\medskip

$iii)$ The space $\cl{U}_0\cap \Delta (\cl{U})$ is a TRO ideal of $\Delta (\cl{U}).$\\
\emph{Proof}\\
Since $\Delta (\cl{U})$ is a TRO $(\cl{U}_0\cap \Delta (\cl{U}))\Delta (\cl{U})^*\Delta (\cl{U})\subset \Delta (\cl{U}).$\\
Using observation $(ii)$ we have that $(\cl{U}_0\cap \Delta (\cl{U}))\Delta (\cl{U})^*\Delta (\cl{U})\subset \cl{U}_0.$\\
It follows that  $(\cl{U}_0\cap \Delta (\cl{U}))\Delta (\cl{U})^*\Delta (\cl{U})\subset \cl{U}_0\cap \Delta (\cl{U}).$\\
Analogously we get $\Delta (\cl{U})\Delta (\cl{U})^*(\cl{U}_0\cap \Delta (\cl{U}))\subset \cl{U}_0\cap \Delta (\cl{U}).$\\
We conclude that the space $\cl{U}_0\cap \Delta (\cl{U})$ is a TRO ideal of $\Delta (\cl{U}).$

\medskip

So there exist projections $Q_i\in \cl{D}_i, i=1,2,$ such that $\cl{U}_0\cap \Delta (\cl{U})=\Delta (\cl{U})Q_1=Q_2\Delta (\cl{U})$ (Remark 2.1).

\medskip

By Theorem 3.1 we have $$\mathcal{U}=\mathcal{U}_0+\Delta(\mathcal{U})=\cl{U}_0+\Delta (\cl{U})Q_1+\Delta (\cl{U})(I-Q_1)=\cl{U}_0+\Delta (\cl{U})(I-Q_1).$$
Clearly $\cl{U}_0\cap \Delta (\cl{U})(I-Q_1)=0.$

\medskip

Similarly one shows that $\cl{U}=\cl{U}_0\oplus (I-Q_2)\Delta (\cl{U})$ and it therefore follows that $\cl{U}=\cl{U}_0\oplus (I-Q_2)\Delta (\cl{U})(I-Q_1). \qquad  \Box$

\begin{remark}The projection $\theta :\cl{U}\rightarrow \cl{U}$ onto $(I-Q_2)\Delta(\mathcal{U})(I-Q_1)$ defined by the decomposition in Theorem 3.3 is a contraction.

\medskip

\emph{Indeed, if $T\in \cl{U},$ as in Theorem 3.1 we have $T=M+S$ where $M\in \Delta (\cl{U}),S\in \cl{U}_0$ and $\|M\|\leq \|T\|$ (see the proof).\\Since $\theta  (T)=(I-Q_2)M(I-Q_1)$, we obtain $\|\theta (T)\|\leq \|T\|.$}
\end{remark}

\bigskip

Let $\mathcal{N}_i=Alg(\mathcal{S}_{i,\phi})=\{T: P^\bot TP=0\;
\text{for all}\; P\in\mathcal{S}_{i,\phi}\},i=1,2,$ and $\mathcal{L}_i=[PTP^\bot:T\in\mathcal{N}_i,P\in\mathcal{S}_{i,\phi}]^{-w^*},\;i=1,2$.

\begin{lemma}i)\;\;$\cl{A}_2\Delta (\cl{U})\cl{A}_1\subset \Delta (\cl{U}).$\\
ii)\;\; $\Delta (\cl{U})^*\cl{A}_2\Delta (\cl{U})\subset \cl{A}_1, \Delta (\cl{U})\cl{A}_1\Delta (\cl{U})^*\subset \cl{A}_2.$\\
iii)\;\;$\cl{U}=\cl{N}_2\cl{U}\cl{N}_1.$\\
iv)\;\;$\cl{U}_0=\cl{N}_2\cl{U}_0\cl{N}_1.$\\
v)\;\;$\cl{U}\cl{L}_1\subset \cl{U}_0,\; \cl{L}_2\cl{U}\subset \cl{U}_0.$\\
vi)\;\;$\Delta (\cl{U})^*\cl{U}\subset \cl{N}_1,\; \cl{U}\Delta (\cl{U})^*\subset \cl{N}_2. $\\
vii)\;\;$\Delta (\cl{U})^*\cl{U}_0\subset \cl{L}_1,\; \cl{U}_0\Delta (\cl{U})^*\subset \cl{L}_2. $
\end{lemma}
\textbf{Proof}

\medskip

Claims $(i),(ii)$ are obvious and $(iii)$ is Lemma 1.1 in \cite{kt}.

\medskip

iv) If $N_1\in \cl{N}_1, N_2\in \cl{N}_2, T\in \cl{U}$ and $P\in \mathcal{S}_{1,\phi }$ then
\begin{align*}N_2\phi (P)TP^\bot N_1=\phi (P)N_2\phi (P)TP^\bot N_1P^\bot \in \cl{U}_0
\end{align*}
since $N_2\phi (P)TP^\bot N_1 \in \cl{U}$ by (iii). Taking the $w^*$ closed linear span we get $\cl{N}_2\cl{U}_0\cl{N}_1\subset \cl{U}_0.$

v) If $N_1\in \cl{N}_1, T\in \cl{U}$ and $P\in \mathcal{S}_{1,\phi } $ then
\begin{align*}TPN_1P^\bot=\phi (P)TPN_1P^\bot\in \cl{U}_0
\end{align*}
  since $TPN_1\in \cl{U}\cl{N}_1\subset \cl{U}$. Taking the $w^*$ closed linear span we get $TK\in \cl{U}_0$ for every $K\in \cl{L}_1.$\\The second inclusion follows by symmetry.

\medskip

vi) Let $M\in \Delta (\cl{U}),T\in \cl{U},P\in\mathcal{S}_{1,\phi }.$ Then $PM^*TP=M^*\phi (P)TP=M^*TP$ so $M^*T\in \cl{N}_1.$\\
Similarly one shows that $TM^*\in \cl{N}_2.$

\medskip

vii) Let $M\in \Delta (\cl{U}),T\in \cl{U},P\in\mathcal{S}_{1,\phi }$ then $M^*\phi (P)TP^\bot=PM^*TP^\bot \in \cl{L}_1$ since $M^*T\in \Delta (\cl{U})^*\cl{U}\subset \cl{N}_1.$Taking the $w^*$ closed linear span we get $M^*S\in \cl{L}_1$ for every $S\in \cl{U}_0.$\\
Similarly one shows that $\cl{U}_0\Delta (\cl{U})^*\subset \cl{L}_2.  \qquad  \Box$

\begin{proposition}
The following are equivalent:
\begin{align*}
i) \qquad & \mathcal{U}  = \mathcal{U}_0.\\
ii) \qquad & \Delta(\cl{U})^*\Delta(\cl{U}) \subset\cl{L}_1\cap \cl{A}_1.\\
iii) \qquad &  \Delta(\cl{U})\Delta(\cl{U})^* \subset\cl{L}_2\cap \cl{A}_2.
\end{align*}
\end{proposition}

\noindent\textbf{Proof }

\medskip

If $\cl{U}=\cl{U}_0$ then $\Delta (\cl{U})\subset \cl{U}_0$, hence $\Delta (\cl{U})^*\Delta (\cl{U})\subset \Delta (\cl{U})^*\cl{U}_0\subset \cl{L}_1$ by the previous lemma.\\
Since $\Delta (\cl{U})^*\Delta (\cl{U})\subset \cl{A}_1$ we get $\Delta (\cl{U})^*\Delta (\cl{U})\subset \cl{L}_1\cap \cl{A}_1$.

\medskip

If conversely $\Delta (\cl{U})^*\Delta (\cl{U})\subset \cl{L}_1\cap \cl{A}_1$, then $\Delta (\cl{U})^*\Delta (\cl{U})(I-Q_1)\subset \cl{L}_1\cap \cl{A}_1,$ so by the previous lemma  $\Delta (\cl{U})\Delta (\cl{U})^*\Delta (\cl{U})(I-Q_1)\subset \cl{U}\cl{L}_1\subset \cl{U}_0.$\\($Q_1$ is the projection in Theorem 3.3).\\
Since $\cl{U}_0\cap \Delta (\cl{U})$ is a TRO ideal of $\Delta (\cl{U})$ (Theorem 3.3) we have that\linebreak $\Delta (\cl{U})\Delta (\cl{U})^*\Delta (\cl{U})Q_1=\Delta (\cl{U})\Delta (\cl{U})^*(\Delta (\cl{U})\cap \cl{U}_0)\subset \Delta (\cl{U})\cap \cl{U}_0\subset \cl{U}_0.$\\
We conclude that $\Delta (\cl{U})\Delta (\cl{U})^*\Delta (\cl{U})\subset \cl{U}_0.$\\
Since $\Delta(\mathcal{U})$ is a TRO its subspace $\Delta(\cl{U})\Delta(\cl{U})^*\Delta(\cl{U})$ is norm-dense \cite{eor}.\\
Therefore $\Delta (\cl{U})\subset \cl{U}_0$ and so $\cl{U}=\cl{U}_0.$

\medskip

The equivalence  $(i)\Leftrightarrow(iii)$ is proved similarly.$\qquad \Box$

\begin{proposition}
The following are equivalent:
\begin{align*}
i) \qquad & \mathcal{U}=\mathcal{U}_0\oplus\Delta(\mathcal{U}). \\
ii) \qquad & \Delta(\mathcal{U})\:(\mathcal{L}_1\cap \cl{A}_1)=0. \\
iii) \qquad & (\mathcal{L}_2\cap \cl{A}_2)\:\Delta(\mathcal{U})=0.
\end{align*}
\end{proposition}

\noindent\textbf{Proof }

\medskip

Note by Lemma 3.5 that $\Delta(\cl{U})(\cl{L}_1\cap \cl{A}_1)\subset \Delta(\cl{U}) \cl{A}_1 \subset \Delta(\cl{U}) $ and \\$\Delta(\cl{U})(\cl{L}_1\cap \cl{A}_1) \subset \cl{U}\cl{L}_1\subset \cl{U}_0.$

\medskip

Thus if the sum $\mathcal{U}=\mathcal{U}_0+\Delta(\mathcal{U})$ is direct then $\Delta(\mathcal{U})\;(\mathcal{L}_1\cap \cl{A}_1)=0.$

\medskip

Suppose conversely that  $\Delta(\mathcal{U})(\mathcal{L}_1\cap \cl{A}_1)=0.$\\
Using again Lemma 3.5 we have that
$(\cl{U}_0\cap \Delta(\cl{U}))^*(\cl{U}_0\cap \Delta(\cl{U}))\subset \Delta (\cl{U})^*\Delta (\cl{U})\\
\subset  \cl{A}_1$ and $(\cl{U}_0\cap \Delta(\cl{U}))^*(\cl{U}_0\cap \Delta(\cl{U}))\subset \Delta (\cl{U})^*\cl{U}_0\subset  \cl{L}_1$ so  \\$(\cl{U}_0\cap \Delta(\cl{U}))^*(\cl{U}_0\cap \Delta(\cl{U}))\subset \cl{L}_1\cap \cl{A}_1$ and so  $(\cl{U}_0\cap \Delta (\cl{U}))(\cl{U}_0\cap \Delta(\cl{U}))^*(\cl{U}_0\cap \Delta(\cl{U}))\\
\subset \Delta (\cl{U})(\cl{L}_1\cap \cl{A}_1)=0.$

But since $\mathcal{U}_0\cap \Delta(\mathcal{U})$ is a TRO (Theorem 3.3), its subspace \linebreak $(\cl{U}_0\cap \Delta(\cl{U}))(\cl{U}_0\cap \Delta(\cl{U}))^*(\cl{U}_0\cap \Delta(\cl{U}))$ is norm-dense \cite{eor}. Therefore $\mathcal{U}_0\cap \Delta(\mathcal{U})=0.$

This shows that (i)and (ii) are equivalent.

\medskip

The proof of the equivalence of (i) and (iii) is analogous. $\qquad \Box$

\section{The diagonal}

Let $\cl{U}, \cl{U}_0, \Delta (\cl{U}), \phi$ be as in section 3 and $\chi =\mathrm{Map} (\Delta (\cl{U})).$

\begin{theorem}\label{38}
There exists a partial isometry $V\in\Delta(\mathcal{U})$ such
that $\Delta(\mathcal{U})=[\cl{A}_2V\cl{A}_1]^{-w^*}$ (recall that
$\cl A_i = (\cl S_{i,\phi})'$).
\end{theorem}
\textbf{Proof}

\medskip

If $T\in \Delta
(\cl{U})$ and $T=U|T|$ is the polar decomposition of $T,$ then $U\in
\Delta (\cl{U})$ and $|T|\in \cl{A}_1:$\;Proposition\;2.6\;in\;\cite{kt}.

\medskip

By Zorn's lemma there exists a maximal family of partial isometries $(V_n)\subset \Delta (\cl{U})$ such that: $V_n^*V_n \bot V_m^*V_m,V_nV_n^* \bot V_mV_m^*$ for $n \ne
m.$

\medskip

Let $V=\sum_{n=1}^{\infty }V_n.$ Then $V$ is a partial isometry in $\Delta (\cl{U}).$

\medskip

First we show that
\begin{equation}\label{3.3}\Delta (\cl{U})=\{T\in B(H_1,H_2): E\in\cl P(\cl{A}_1')
,\, F\in\cl P( \cl{A}_2'),\,  FVE=0\; \Rightarrow\;FTE=0 \}\end{equation}

\medskip

Let T be such that, if $FVE=0$ for $E\in \cl{P}(\cl{A}_1^\prime)$ and $F\in \cl{P}(\cl{A}_2^\prime),$ then $FTE=0.$
Since $\phi (P)^\bot VP=\phi (P)VP^\bot=0$ for every $P\in
\cl{S}_{1,\phi}$ and $\cl{S}_{i,\phi}\subset \cl{A}_i^\prime,
i=1,2$ we have $\phi (P)^\bot TP=\phi (P)TP^\bot=0$ for every $P\in
\cl{S}_{1,\phi}$  so $T\in \Delta (\cl{U}).$

\medskip

For the converse let $T\in \Delta (\cl{U})$ and $T=U|T|$ be the polar decomposition of T.

\medskip

If $E\in \cl{P}(\cl{A}_1^\prime)$,$F\in \cl{P}(\cl{A}_2^\prime)$ are such that $FVE=0$, since $|T|\in \cl{A}_1,$ we have $FTE=FU|T|E=FUE|T|.$

\medskip

Hence it suffices to show that $FUE=0.$

\medskip

We observe that:
\begin{align*}
V^*V(FUE)^*FUE &=(V^*V)EU^*FUE=E(V^*V)U^*FUE \quad (V^*V\in \cl A_1)\\
&=EV^*(VU^*)FUE =EV^*F(VU^*)UE \quad (VU^*\in \cl A_2)\\
&=0VU^*UE=0 \end{align*} hence
\begin{equation}\label{3.4}(FUE)^*FUE \leq I-V^*V. \end{equation}
Similarly, one shows that

\begin{equation}\label{3.5} FUE(FUE)^*\leq I-VV^*. \end{equation}

\medskip

Since $FUE$ is a partial isometry in $\Delta (\cl{U}),$ the maximality of $V$ and (4.2),(4.3) imply that $FUE=0.$

\medskip

Thus claim (4.1) holds.

\medskip

Let $\cl{M}=[\cl{A}_2V\cl{A}_1]^{-w^*}.$

\medskip

We observe that $\cl{M}$ is a TRO which is contained in $\Delta (\cl{U}).$ Since $\cl{M}$ is $w^*$ closed, it is reflexive.

\medskip

If $\zeta =\mathrm{Map} (\cl{M})$ then for every projection $P$,$$ \zeta  (P)=[A_2VA_1Py :A_i\in \cl{A}_i, i=1,2,y\in H_1]^{-}.$$  We observe that $\zeta  (P)\in \cl{A}_2^\prime$ for every projection P so $\cl{S}_{2,\zeta  }\subset \cl{A}_2^\prime.$ Similarly  if $\zeta ^* =\mathrm{Map} (\cl{M}^*)$\;then\; $\cl{S}_{2,\zeta ^* }\subset \cl{A}_1^\prime$ but $\cl{S}_{1,\zeta  }=\{P^\bot :P\in \cl{S}_{2,\zeta ^* }\}$ so we have that $\cl{S}_{1,\zeta }\subset \cl{A}_1^\prime.$

\medskip

Now since $V\in \cl{M}$ we conclude that $\zeta  (P)^{\bot}VP=0$ for every   $P\in\mathcal{S}_{1,\zeta  }.$

\medskip

 From claim (4.1) we obtain $\zeta  (P)^{\bot}\Delta (\cl{U})P=0 $ for every $P\in\mathcal{S}_{1,\zeta  },$  so since $\cl{M}$ is reflexive $\Delta (\cl{U})\subset \cl{M}. \qquad \Box$

\medskip

By the previous theorem it follows that if $\cl{M}$ is a $w^*-$closed TRO masa bimodule and $\zeta =\mathrm{Map} (\cl{M})$ then there exists a partial isometry $V\in \cl{M}$ so that $\cl{M}=[(\cl{S}_{2,\zeta })^\prime V (\cl{S}_{1,\zeta })^\prime]^{-w^*}.$\\
But we shall prove a stronger result:

\begin{theorem}Let $\cl{M}$ a $w^*-$closed TRO masa bimodule and $\cl{B}_1=[\cl{M}^*\cl{M}]^{-w^*}\\\cl{B}_2=[\cl{M}\cl{M}^*]^{-w^*}.$ Then there exists a partial isometry $V$ such that $\cl{M}=[\cl{B}_2V\cl{B}_1]^{-w^*}.$
\end{theorem}
\textbf{Proof}

\medskip

Let $\cl{D}_i\subset B(H_i), i=1,2$ be masas such that $\cl{D}_2\cl{M}\cl{D}_1\subset \cl{M}$ and put $\zeta =\mathrm{Map} (\cl{M}).$\\
We shall prove that $\cl{B}_2^\prime \cl{M}\cl{B}_1^\prime \subset \cl{M}.$\\
In \cite{kt}, Theorem 2.10 it is shown that $$\cl{B}_2^\prime=(\cl{M}\cl{M}^*)^\prime\subset \cl{D}_2|_{\zeta (I)}\oplus B(\zeta (I)^\bot (H_2))$$ and $$\cl{B}_1^\prime=(\cl{M}^*\cl{M})^\prime\subset \cl{D}_1|_{\zeta^* (I)}\oplus B(\zeta^* (I)^\bot (H_1)).$$
So it suffices to show that $$(\cl{D}_2|_{\zeta (I)}\oplus B(\zeta (I)^\bot (H_2)))\;\cl{M}\;(\cl{D}_1|_{\zeta^* (I)}\oplus B(\zeta^* (I)^\bot (H_1)))\subset \cl{M}.$$
But this is true because $\cl{D}_2\cl{M}\cl{D}_1\subset \cl{M}, \zeta (I)\in \cl{D}_2, \zeta^* (I)\in \cl{D}_1$ and $\cl{M}=\zeta (I)\cl{M}\zeta^* (I).$

\medskip

Now, we shall follow the proof of the previous theorem:\\
By Zorn's lemma there exists a maximal family of partial isometries $(V_n)\subset \cl{M}$ such that: $V_n^*V_n \bot V_m^*V_m,V_nV_n^* \bot V_mV_m^*$ for $n \ne
m.$\\
Let $V=\sum_{n=1}^{\infty }V_n.$ Then $V$ is a partial isometry in $\cl{M}.$\\
We shall show that
\begin{equation}\label{3.4}\cl{M}\subset \{T\in B(H_1,H_2): E\in\cl P(\cl{B}_1')
,\, F\in\cl P( \cl{B}_2'),\,  FVE=0\; \Rightarrow\;FTE=0 \}\end{equation}

Let $T\in \cl{M}$ and $T=U|T|$ be the polar decomposition of $T.$ Then $|T|\in (\cl{M}^* \cl{M})^{\prime \prime}$ and $U \in \cl{M},$ (Proposition 2.6 in \cite{kt}).\\
If $E\in\cl P(\cl{B}_1')
,\, F\in\cl P( \cl{B}_2')$ are such that $FVE=0,$ since $|T|\in (\cl{M}^* \cl{M})^{\prime \prime}$ and $E\in \cl{B}_1^\prime=(\cl{M}^* \cl{M})^{\prime },$ we have $FTE=FU|T|E=FUE|T|.$ Hence it suffices to show that $FUE=0.$

\medskip

As in the proof of the previous theorem we have that $V^*V\bot (FUE)^*(FUE)$ and $VV^*\bot (FUE)(FUE)^*.$\\
But $FUE\in \cl{B}_2^\prime \cl{M} \cl{B}_1^\prime \subset \cl{M},$ so by the maximality of $V$ we have that $FUE=0.$

\medskip

Let $\cl{W}=[\cl{B}_2 V \cl{B}_1]^{-w^*}.$ We observe that $\cl{W}\subset \cl{M}.$\\
For the converse, we follow the proof of the previous theorem and we use the relation (4.4) $\qquad  \Box$

An alternative proof of the previous theorem was communicated to us by I. Todorov, based on his paper \cite{todo}.

\begin{theorem}\label{43} The semilattices of $\Delta (\cl{U})$ are the following:
$$\cl{S}_{1,\chi }=\chi ^*(I)^\bot \oplus \chi ^*(I)\cl{P}((\cl{S}_{1,\phi })^{\prime \prime})$$
$$\cl{S}_{2,\chi }=\chi (I)\cl{P}((\cl{S}_{2,\phi })^{\prime \prime}).$$
The map $\chi :\cl{S}_{1,\chi }\longrightarrow \cl{S}_{2,\chi }$
is such that
\begin{equation}\label{xi}
\chi (\chi ^*(I)^\bot \oplus \chi ^*(I)Q)=\chi (I)\phi (Q) \quad
\text{for every} \; Q\in \cl{S}_{1,\phi }.
\end{equation}
\end{theorem}
\textbf{Proof}

\medskip

i) In Theorem \ref{38} we showed that there exists a partial
isometry $V$ in $\Delta (\cl{U})$ such that $\Delta
(\cl{U})=[(\cl{S}_{2,\phi })^{\prime}V(\cl{S}_{1,\phi
})^{\prime}]^{-w^*}$.

\medskip

So if $P\in \cl{S}_{1,\chi }$ then $\chi (P)$ is the projection onto  $[(\cl{S}_{2,\phi })^{\prime}V(\cl{S}_{1,\phi })^{\prime}P(H_1)]^{-}$.

\medskip

We conclude that $\chi (P)\in (\cl{S}_{2,\phi })^{\prime \prime}.$ Hence $\cl{S}_{2,\chi }\subset (\cl{S}_{2,\phi })^{\prime \prime}.$

\medskip

If $H$ is a Hilbert space, $\cl{B}$ is a subset of $B(H)$ and $Q$ a projection in  $\cl{B}^\prime$ the set $\{T|_{Q(H)}:T\in \cl{B}\}$ is denoted by $\cl{B}|_Q$.

\medskip

We have shown that $(\cl{S}_{2,\chi})^{\prime \prime}|_{\chi (I)}\subset (\cl{S}_{2,\phi })^{\prime \prime}|_{\chi (I)}.$

\medskip

Let $P\in \cl{S}_{1,\phi }$ then $\Delta (\cl{U})P=\phi (P)\Delta (\cl{U}).$ Hence $\chi (P)=\phi (P)\chi (I).$

\medskip

So $\chi (I)\cl{S}_{2,\phi }\subset \cl{S}_{2,\chi }$ hence, $(\cl{S}_{2,\phi })^{\prime \prime}|_{\chi (I)}\subset (\cl{S}_{2,\chi })^{\prime \prime}|_{\chi (I)}.$

\medskip

We proved that $$(\cl{S}_{2,\phi })^{\prime \prime}|_{\chi (I)} = (\cl{S}_{2,\chi })^{\prime \prime}|_{\chi (I)}.$$

\medskip

Since $\Delta (\cl{U})$ is a TRO, using Theorem 2.10 in \cite{kt} (see the proof) we have that $$\cl{S}_{2,\chi }|_{\chi(I)} =\cl{P}((\cl{S}_{2,\chi })^{\prime \prime}|_{\chi (I)}).$$

\medskip

It follows that $$\cl{S}_{2,\chi }=\chi (I)\cl{P}((\cl{S}_{2,\phi })^{\prime \prime}).$$

\medskip

Applying this to $\Delta (\cl{U})^*=\Delta (\cl{U}^*),$
$$\cl{S}_{2,\chi^* }=\chi^* (I)\cl{P}((\cl{S}_{2,\phi^* })^{\prime \prime}).$$

\medskip

Since $\cl{S}_{1,\phi}=\{Q^\bot: Q\in \cl{S}_{2,\phi^*}\},$ see the introduction, we have that $$\cl{S}_{2,\chi^* }=\chi^* (I)\cl{P}((\cl{S}_{1,\phi })^{\prime \prime}).$$

\medskip

But \begin{align*}\cl{S}_{1,\chi}=\{Q^\bot: Q\in \cl{S}_{2,\chi^*}\}&=\{(\chi ^*(I)Q)^\bot: Q\in \cl{P}((\cl{S}_{1,\phi })^{\prime \prime})\}\\&=\{\chi ^*(I)^\bot \oplus \chi ^*(I)Q: Q\in \cl{P}((\cl{S}_{1,\phi })^{\prime \prime})\}.\end{align*}

\medskip

ii) If $Q\in \cl{S}_{1,\phi }$ then

\begin{align*}
\chi (\chi ^*(I)^\bot \oplus \chi ^*(I)Q)&=\chi (\chi ^*(I)Q)\qquad (\chi (\chi ^*(I)^\bot )=0)\\
&=\chi (Q)\quad \quad \quad \quad(\Delta (\cl{U})\chi ^*(I)=\Delta (\cl{U}))\\
&=\phi (Q)\chi (I).\qquad (\Delta (\cl{U})Q=\phi (Q)\Delta (\cl{U}))\end{align*}
$\qquad \Box$

\begin{remark}The smallest ortholattice containing the commutative family $\chi (I)\cl{S}_{2,\phi}$ is easily seen to be
$\chi (I)\cl{P}((\cl{S}_{2,\phi})'')$, which equals
$\cl{S}_{2,\chi}$; similarly the family $\chi ^*(I)^\bot \oplus
\chi ^*(I)\cl{S}_{1,\phi}$ generates the complete ortho-lattice
$\cl{S}_{1,\chi }$.

\medskip

Therefore, since $\chi |_{\cl{S}_{1,\chi }}$ is a complete
ortho-lattice isomorphism (Theorem 2.10 in \cite{kt}) equality
 (\ref{xi}) determines the map $\chi.$

\end{remark}

\begin{proposition}The families $\chi ^*(I)\cl{S}_{1,\phi }$ and $\chi (I)\cl{S}_{2,\phi }$ are complete lattices and the map $$\vartheta :\chi ^*(I)\cl{S}_{1,\phi }\rightarrow  \chi (I)\cl{S}_{2,\phi }: \vartheta (\chi ^*(I)P)=\chi (I)\phi (P)$$ is a complete lattice isomorphism.
\end{proposition}
\textbf{Proof}

\medskip

We use Theorem \ref{43} and the fact \cite{kt} that the map $\chi
|_{\cl{S}_{1,\chi }}$ is a complete ortholattice isomorphism .

\medskip

Let $(P_i)_{i\in I}\subset \cl{S}_{1,\phi }$. We claim that
\begin{equation}\label{sx45}
\wedge _{i\in I}\chi (I)\phi (P_i)=\chi (I)\phi (\wedge _{i\in
I}P_i).
\end{equation}
Indeed, by (\ref{xi}),
\begin{align*}
\wedge _{i\in I}\chi (I)\phi (P_i)= & \wedge _{i\in I}\chi (\chi
^*(I)^\bot \oplus \chi ^*(I)P_i) \\
= & \chi (\wedge _{i\in I}(\chi ^*(I)^\bot \oplus \chi ^*(I)P_i))=
\chi (\chi ^*(I)^\bot \oplus \chi ^*(I)(\wedge _{i\in
I}P_i)).
\end{align*}
Since $\wedge _{i\in I}P_i \in \cl{S}_{1,\phi }$ we get that $\chi
(\chi ^*(I)^\bot \oplus \chi ^*(I)(\wedge _{i\in I}P_i))=\chi
(I)\phi (\wedge _{i\in I}P_i)$ again using (\ref{xi}).

\medskip

By (\ref{sx}), there exist $(Q_i)_{i\in I}\subset \cl{S}_{1,\phi^* }$ such that $\phi ^*(Q_i)^\bot =P_i$ for every $i\in I.$\\
We shall prove that
\begin{equation}\label{sx46}
\vee _{i\in I}\chi ^*(I)P_i=\chi ^*(I)(\phi ^*(\wedge _{i\in
I}Q_i))^\bot.
\end{equation}
 Since $\Delta (\cl{U}^*)=\Delta
(\cl{U})^*$ we have that $\chi ^*=Map(\Delta (\cl{U}^*))$ and so
applying  equation (\ref{sx45}) to $\chi^*$ we have that
 \begin{align*}
 \wedge _{i\in I}\chi^* (I)\phi^* (Q_i)= & \chi^*
(I)\phi^* (\wedge _{i\in I}Q_i) \; \Rightarrow \\
\vee  _{i\in I}(\chi^* (I)\phi^* (Q_i))^\bot = & (\chi^* (I)\phi^*
(\wedge _{i\in I}Q_i)^\bot \; \Rightarrow \\
\vee _{i\in I}(\chi^* (I)^\bot \oplus \chi ^*(I)(\phi^* (Q_i))^\bot)
= & \chi^* (I)^\bot \oplus \chi ^*(I)(\phi^* (\wedge _{i\in
I}Q_i))^\bot \; \Rightarrow \\
\vee _{i\in I}(\chi ^*(I)^\bot \oplus \chi ^*(I)P_i)= & \chi^*
(I)^\bot \oplus \chi ^*(I)(\phi^* (\wedge _{i\in I}Q_i))^\bot \;
\Rightarrow \\
\vee _{i\in I}\chi ^*(I)P_i= & \chi ^*(I)(\phi^* (\wedge _{i\in
I}Q_i))^\bot.
 \end{align*}

From equalities $(4.6)$ and $(4.7)$ we conclude that the families
$\chi ^*(I)\cl{S}_{1,\phi },\\ \chi (I)\cl{S}_{2,\phi }$ are
complete lattices.

\medskip

Since $\chi (\chi ^*(I)^\bot \oplus \chi ^*(I)Q)=\chi (I)\phi (Q)$ for every $Q\in \cl{S}_{1,\phi }$
and $\chi|_{\cl{S}_{1,\chi }} $ is $1-1$ the map $\vartheta$ is a bijection.\\
It remains to show that $\vartheta $ is sup and inf continuous.

\medskip

Let   $(P_i)_{i\in I}\subset \cl{S}_{1,\phi }$ and $(Q_i)_{i\in I}\subset \cl{S}_{1,\phi^* }$ be
such that $\phi ^*(Q_i)^\bot =P_i$, equivalently by equation (1.1) $\phi (P_i)^\bot =Q_i$ for every $i\in I.$\\
Then, since $\wedge _{i\in I}P_i\in \cl{S}_{1,\phi }$, by the
definition of $\vartheta$ we have
\begin{align*}
\vartheta (\wedge _{i\in I}\chi ^*(I)P_i)= &\ \vartheta (\chi
^*(I)(\wedge _{i\in I}P_i))=\chi (I)\phi (\wedge _{i\in I}P_i)\\
=&\ \wedge _{i\in I}\chi (I)\phi (P_i)=\wedge _{i\in I}\vartheta
(\chi ^*(I)P_i).
\end{align*}
Using equations $(4.7)$ and $(1.1)$ we have that
\begin{align*}
\vartheta (\vee  _{i\in I}\chi ^*(I)P_i)= & \ \vartheta (\chi
^*(I)(\phi ^*(\wedge _{i\in I}Q_i))^\bot )
= \chi (I)\phi ((\phi ^*(\wedge _{i\in I}Q_i))^\bot )\\
=& \ \chi (I)(\wedge _{i\in I}Q_i)^\bot =\vee _{i\in I}\chi
(I)Q_i^\bot \\
=& \ \vee _{i\in I}\chi (I)\phi (P_i)=\vee  _{i\in I}\vartheta
(\chi ^*(I)P_i). \qquad \Box
\end{align*}

\section{The space $\cl{U}_0$ is reflexive.}

Let $\cl{U}, \cl{U}_0, \Delta (\cl{U}), \phi$ be as in section 3 and $\chi =\mathrm{Map} (\Delta (\cl{U})), \psi =\mathrm{Map} (\cl{U}_0).$

\begin{lemma}If $\Delta (\cl{U})$ is essential, i.e. $\chi (I)=I, \chi ^*(I)=I,$ then $\cl{S}_{1,\psi }\subset \cl{S}_{1,\phi}$ and $\cl{S}_{2,\psi }\subset \cl{S}_{2,\phi}.$
\end{lemma}
\textbf{Proof}

\medskip

Since $\chi (I)=I$ we have $\phi (I)=I,$ so by Proposition 4.5, $\cl{S}_{2,\phi }$ is a C.S.L.\\
Since $\chi^* (I)=I,$ $\cl{S}_{2,\phi^* }$ is a C.S.L. and so $\cl{S}_{1,\phi }$ is a C.S.L.

\medskip

If $E$ is a projection, then $Alg(\cl{S}_{2,\phi })\cl{U}_0E\subset \cl{U}_0E$ (Lemma 3.5).\\
It follows that $\psi (E)^\bot Alg(\cl{S}_{2,\phi })\psi (E)=0.$ Hence $\psi (E)\in Lat(Alg(\cl{S}_{2,\phi })).$ Since commutative subspace lattices are reflexive \cite{arv}, it follows that $\psi (E)\in \cl{S}_{2,\phi }.$ We get that $\cl{S}_{2,\psi }\subset \cl{S}_{2,\phi}.$ \\
Analogously $\cl{U}_0Alg(\cl{S}_{1,\phi })\subset Alg(\cl{S}_{1,\phi })$ so $Alg(\cl{S}_{1,\phi }^\bot ) \cl{U}_0^*\subset \cl{U}_0^*.$ As above we obtain $\cl{S}_{2,\psi^* }\subset \cl{S}_{1,\phi}^\bot$ hence $\cl{S}_{1,\psi}\subset \cl{S}_{1,\phi}. \qquad \Box$

\begin{theorem}The space $\cl{U}_0$ is reflexive
\end{theorem}
\textbf{Proof}

\medskip

Firstly, we suppose that $\Delta (\cl{U})$ is essential  ($\chi (I)=I, \chi ^*(I)=I$).\\
Now, by Theorem 4.3 we have that $\cl{S}_{1,\chi }=\cl{P}((\cl{S}_{1,\phi })^{\prime \prime}), \cl{S}_{2,\chi }=\cl{P}((\cl{S}_{2,\phi })^{\prime \prime})$ and $\chi |_{\cl{S}_{1,\phi }}=\phi .$

\medskip

If $E\in \cl{S}_{1,\phi },$ then $\phi (E), \psi (E) \in \cl{P}((\cl{S}_{2,\phi })^{\prime \prime})$ so there exists a unique $F\in \cl{P}((\cl{S}_{1,\phi })^{\prime \prime})$ such that $\chi (F)=\phi (E)-\psi (E).$\\
We observe that $\chi (F)\leq \phi (E)=\chi (E).$ Since $\chi $ is a lattice isomorphism $F\leq E$ and so $\psi (F)\leq \psi (E);$ therefore $\chi (F)\bot \psi (F).$

Since  $\chi =\mathrm{Map} (\Delta (\cl{U}))$ and $\psi =\mathrm{Map} (\cl{U}_0)$ we obtain that\\ $\Delta (\cl{U})F(H_1)\bot \mathrm{Ref}(\cl{U}_0)F(H_1)$ and so $\Delta (\cl{U})F \cap  \mathrm{Ref}(\cl{U}_0)F=0.$\\
By Theorem 3.1 $\cl{U}=\cl{U}_0+\Delta (\cl{U}),$ hence $\cl{U}F=\mathrm{Ref}(\cl{U}_0)F\oplus \Delta (\cl{U})F$ and $\cl{U}F=\cl{U}_0 F\oplus \Delta (\cl{U})F.$\\
It follows that $\cl{U}_0 F=\mathrm{Ref}(\cl{U}_0)F$ and so $\cl{U}_0F$ is reflexive.

\medskip

Let $$P=\vee \{F\in \cl{P}((\cl{S}_{1,\phi })^{\prime \prime}): \chi (F)=\phi (E)-\psi (E), E\in  \cl{S}_{1,\phi }\}.$$
By the previous arguments the space $\cl{U}_0P$ is reflexive.\\
Since $\chi $ is $\vee -$continuous we have that $$\chi (P)=\vee \{\phi (E)-\psi (E), E\in  \cl{S}_{1,\phi }\}.$$
Let $Q=\chi (P)^\bot $ then $Q\phi (E)=Q\psi (E)$ for all $E\in \cl{S}_{1,\phi }.$ Therefore, it follows that
$$Q\cl{U}=\{T: Q\phi (E)^\bot TE=0\;\text{for\;all}\;E\in \cl{S}_{1,\phi }\}=$$
 $$=\{T: Q\psi  (E)^\bot TE=0\;\text{for\;all}\;E\in \cl{S}_{1,\phi }\}.$$
Using the previous lemma ($\cl{S}_{1,\psi }\subset \cl{S}_{1,\phi}$) we obtain that $Q\cl{U}$ is contained in the space:$$\{T: Q\psi  (E)^\bot TE=0\;\text{for\;all}\;E\in \cl{S}_{1,\psi  }\}=$$
$$=Q=\mathrm{Ref}(\cl{U}_0)=\mathrm{Ref}(Q\cl{U}_0)\subset Q\cl{U}.$$
We proved that $Q\cl{U}=\mathrm{Ref}(Q\cl{U}_0).$

\medskip

Katavolos and Todorov \cite{kt} have proved that $\Delta (\cl{U})\subset (\cl{U})_{min}$ where $(\cl{U})_{min}$ is the smallest $w^*-$closed masa bimodule such that $\mathrm{Ref}((\cl{U})_{min})=\cl{U}.$\\
So $Q\Delta (\cl{U})\subset Q(\cl{U})_{min}=(Q\cl{U})_{min}.$ But since $Q\cl{U}_0$ is a $w^*-$closed masa bimodule such that $\mathrm{Ref}(Q\cl{U}_0)=Q\cl{U}$ it follows that $Q\Delta (\cl{U})\subset Q\cl{U}_0.$\\
Now $Q\Delta (\cl{U})=\chi (P)^\bot \Delta (\cl{U})=\Delta (\cl{U})P^\bot,$ hence $\Delta (\cl{U})P^\bot\subset \cl{U}_0.$\\
So $\cl{U}=\cl{U}_0+\Delta (\cl{U})P^\bot+\Delta (\cl{U})P=\cl{U}_0+\Delta (\cl{U})P$ and so $\cl{U}P^\bot =\cl{U}_0P^\bot .$\\
We conclude that $\cl{U}_0P^\bot $ is reflexive. Since  $\cl{U}_0P$ is reflexive too,\;\;$\cl{U}_0$ is reflexive.

\medskip

Now, relax the assumption that $\Delta (\cl{U})$ is essential.

Let $\cl{W}=\chi (I)\cl{U}|_{\chi ^*(I)}.$ This is a masa bimodule in $B(\chi ^*(I)(H_1),\chi (I)(H_2)).$\\
We have that $$\cl{W}=\{T: \chi (I)\phi (L)^\bot TL|_{\chi ^*(I)}=0\;\;\text{for\;all}\;L\in \cl{S}_{1,\phi }\}.$$
By Proposition 4.5 the families $\cl{S}_{1,\phi }|_{\chi ^*(I)}, \cl{S}_{2,\phi }|_{\chi (I)}$ are complete lattices and the map $\cl{S}_{1,\phi }|_{\chi ^*(I)}\rightarrow  \cl{S}_{2,\phi }|_{\chi (I)}: P|_{\chi ^*(I)}\rightarrow \phi (P)|_{\chi (I)}$ is a complete lattice isomorphism.\\
By the Lifting theorem of J.Erdos \cite{erd} it follows that the (semi)lattices of $\cl{W}$ are the families $\cl{S}_{1,\phi }|_{\chi ^*(I)}, \cl{S}_{2,\phi }|_{\chi (I)}.$\\
Therefore, $\cl{W}_0=[\chi (I)\phi (L)TL^\bot |_{\chi ^*(I)}:\; T\in \cl{W}, L\in \cl{S}_{1,\phi }]^{-w^*}=\chi (I)\cl{U}_0|_{\chi ^*(I)}.$

\medskip

By the proof in the essential case we have that the spase $\chi (I)\cl{U}_0 \chi ^*(I)$ is reflexive.\\
But $\chi (I)^\bot \cl{U}=\chi (I)^\bot \cl{U}_0$ and $\cl{U}\chi ^*(I)^\bot =\cl{U}_0\chi ^*(I)^\bot $ so the spaces
$\chi (I)^\bot \cl{U}_0$ and $\cl{U}_0\chi ^*(I)^\bot $ are reflexive.\\
Finally the space  $\cl{U}_0$ is reflexive. $\qquad \Box$

\bigskip

For the rest of this section let $\cl{S}$ be a C.S.L.\;\;\; $\cl{U}=Alg(\cl{S}), \cl{J}=[PTP^\bot : T\in \cl{U}, P\in \cl{S}]^{-\|\cdot\|}, Rad(\cl{U})$ be the radical of $\cl{U},$ $\cl{U}_0=\cl{J}^{-w^*}, \psi =\mathrm{Map} (\cl{U}_0).$
It is known that $\cl{J}\subset Rad(\cl{U}).$ The equality $\cl{J}= Rad(\cl{U})$ is an open problem (Hopenwasser's conjecture), \cite{hop}, \cite{do}.\\
I.Todorov \cite{tod} has proved that $\cl{J}$ and $Rad(\cl{U})$ have the same reflexive hull.\\We improve this by showing the next corollary.

\begin{corollary}The spaces $\cl{J}$ and $Rad(\cl{U})$ have the same $w^*-$closure.
\end{corollary}
\textbf{Proof}

\medskip

$$\cl{U}_0=\cl{J}^{-w^*}\subset Rad(\cl{U})^{-w^*}\subset \mathrm{Ref}(Rad(\cl{U}))=\mathrm{Ref}(\cl{J})= \cl{U}_0.$$ $\qquad \Box$

\begin{corollary}$Rad(\cl{U})^{-w^*}=\{T: \psi (E)^\bot TE=0 \;\text{for\;all}\;E\in \cl{S}\}.$
\end{corollary}
\textbf{Proof}

\medskip

$Rad(\cl{U})^{-w^*}=\cl{U}_0=\{T: \psi  (E)^\bot TE=0\;\text{for\;every\;projection}\;E\}\subset\\ \{T: \psi  (E)^\bot TE=0\;\text{for\;all}\;E\in \cl{S}\}.$\\
Using Lemma 5.1 the last space is contained in the space:\\
$\{T: \psi  (E)^\bot TE=0\;\text{for\;all}\;E\in \cl{S}_{1,\psi }\}=\cl{U}_0=Rad(\cl{U})^{-w^*}. \qquad \Box$

\medskip

Now we are ready to give the form of the decomposition of $\cl{U}$ in the case that $\cl{U}$ is a C.S.L. algebra:

\begin{proposition}Let $Q=\vee \{E-\psi (E): E\in \cl{S}\}$ then $$\cl{U}=Rad(\cl{U})^{-w^*}\oplus Q\cl{S}^{\prime}.$$
\end{proposition}
\textbf{Proof}

\medskip

We observe that  $Q^\bot E=Q^\bot \psi (E)$ for all $E\in \cl{S},$ so we have:
$$Q^\bot \cl{U}=\{T: Q^\bot E^\bot TE=0\;\text{for\;all}\;E\in \cl{S}\}$$
 $$=\{T: Q^\bot \psi  (E)^\bot TE=0\;\text{for\;all}\;E\in \cl{S}\}.$$
By the previous corollary the last space is the space $Q^\bot Rad(\cl{U})^{-w^*}.$\\
So we have that $Q^\bot \cl{S}^{\prime}\subset Q^\bot Rad(\cl{U})^{-w^*}\subset Rad(\cl{U})^{-w^*}.$\\
Since $\cl{U}=Rad(\cl{U})^{-w^*}+\cl{S}^{\prime}$ we have $\cl{U}=Rad(\cl{U})^{-w^*}+ Q\cl{S}^{\prime}.$

\medskip

It suffices to show that $Rad(\cl{U})^{-w^*}\cap Q\cl{S}^{\prime}=0.$\\
Let $E\in \cl{S}$ and $T\in \cl{U}_0\cap (E-\psi (E))\cl{S}^{\prime}$ then $T=(E-\psi (E))T\\=\psi (E)^\bot ET=\psi (E)^\bot TE=0,$ because $T\in \cl{U}_0.$\\
If $T\in \cl{U}_0\cap Q\cl{S}^{\prime}\Rightarrow (E-\psi (E))T\in \cl{U}_0\cap (E-\psi (E))\cl{S}^{\prime}=0.$\\
So $(E-\psi (E))T=0$ for all $E\in \cl{S}.$\\
But $T=(\vee \{E-\psi (E): E\in \cl{S}\})T.$ It follows that $T=0 \qquad \Box$

\section{Decomposition of compact operators in reflexive masa bimodules}
Let $\cl{U},\;\cl{U}_0,\;\Delta (\cl{U}),\;\phi ,\;\cl{D}_1,\;\cl{D}_2,\;Q_1$ be as in section 3 and $\chi=\mathrm{Map} (\Delta (\cl{U}))$ .

\bigskip

We denote by $\cl{K}$ the set of compact operators and by $C_p$ the set of p-Schatten class
 operators in $B(H_1,H_2).$

\begin{proposition}If\; $T\in R_1(\mathcal{U})$, there exist $L\in R_1(\Delta(\mathcal{U}))$ and \\ $S\in[R_1(\mathcal{U}_0)]^{-\|\;\|_1}$ such that $T=L+S.$
\end{proposition}
\textbf{Proof}

\medskip

Write $\cl{U}=\{X : \phi (P_n)^\bot XP_n=0\; \text{for all}\; n\in \mathbb{N}\}$ for an approriate sequence $(P_n)\subset \cl{S}_{1,\phi }$ and let $T\in R_1(\mathcal{U}).$

\medskip

As in the proof of Theorem 3.1
$$T=L_1+\phi (P_1)TP_1^\bot,\text{where}\; L_1=\phi (P_1)TP_1+\phi (P_1)^\bot TP_1^\bot.$$
Since $\phi (P_1)^\bot TP_1=0$ and $T$ has rank 1 either $\phi (P_1)^\bot T=0$ or $TP_1=0$, hence either $L_1=\phi (P_1)TP_1$ or $L_1=\phi (P_1)^\bot TP_1^\bot.$
$$L_1=L_2+\phi (P_2)L_1P_2^\bot,\text{where}\; L_2=\phi (P_2)L_1P_2+\phi (P_2)^\bot L_1P_2^\bot.$$
Since $\phi (P_2)^\bot L_1P_2=0,$ either $L_2=\phi (P_2)L_1P_2$ or $L_2=\phi (P_2)^\bot L_1P_2^\bot.$
Similarly $$L_{n-1}=L_n+\phi (P_n)L_{n-1}P_n^\bot,\text{where}\; L_n=\phi (P_n)L_{n-1}P_n+\phi (P_n)^\bot L_{n-1}P_n^\bot .$$As before, either $L_n=\phi (P_n)L_{n-1}P_n$ or $L_n=\phi (P_n)^\bot L_{n-1}P_n^\bot$ for all $n\in \mathbb{N}.$

\medskip

We conclude that there exist projections $(Q_n)\subset \cl{D}_2, (R_n)\subset \cl{D}_1$ such that $L_n= (\wedge _{i=1}^n Q_i) T (\wedge _{i=1}^n R_i) ,n\in \mathbb{N}.$

\medskip

We observe that $T=L_n+M_n$ where $M_n=\phi (P_1)TP_1^\bot+\phi (P_2)L_1P_2^\bot+...+\phi (P_n)L_{n-1}P_n^\bot,$ $ n\in \mathbb{N}.$

\medskip

 Since $\wedge _{i=1}^n Q_i \stackrel {sot}{\rightarrow}\wedge _{i=1}^{\infty } Q_i,$\;$\wedge _{i=1}^n R_i \stackrel {sot}{\rightarrow}\wedge _{i=1}^{\infty } R_i$ and T has rank 1  $$L_n \stackrel{\|\;\|_1}{\rightarrow }(\wedge _{i=1}^{\infty } Q_i) T (\wedge _{i=1}^{\infty } R_i)=L, \text{say.} $$

\medskip

Now $\phi (P_i)^\bot L_nP_i=\phi (P_i)L_nP_i^\bot=0, i=1,2,...n$ for all $n\in \mathbb{N},$ therefore $\phi (P_i)^\bot LP_i=\phi (P_i)LP_i^\bot=0$ for all $i\in \mathbb{N}.$

\medskip

Thus $L\in R_1(\Delta (\cl{U})).$

\medskip

We have $M_n=T-L_n \stackrel{\|\;\|_1}{\rightarrow }T-L=S\in [R_1(\cl{U}_0)]^{-\|\;\|_1}. \qquad \Box$

\begin{proposition}$\cl{U}_0 \subset (R_1(\Delta (\cl{U}))^*)^0.$
\end{proposition}
\textbf{Proof}

\medskip

Let $T\in\cl{U},P\in\cl{S}_{1,\phi},R\in R_1(\Delta (\cl{U})).$ Then $$tr(\phi (P)TP^\bot R^*)=\linebreak tr(T(\phi (P)RP^\bot)^*)=tr(T0)=0.$$

\medskip

Taking the $w^*$ closed linear span we get $tr(SR^*)=0$ for every $S\in \cl{U}_0.\qquad \Box$

\begin{proposition}\label{42}$i)\;R_1(\Delta (\cl{U}))\subset \Delta (\mathcal{U})(I-Q_1).$\\
$ii)\;\Delta (\cl{U})\cap \cl{K}=[R_1(\Delta (\cl{U}))]^{-\|\cdot\|}\subset \Delta (\mathcal{U})(I-Q_1).$
\end{proposition}
\textbf{Proof}

\medskip

Let $R\in R_1(\Delta (\cl{U}))$ then as in Theorem 3.3 $RQ_1\in \Delta (\cl{U})Q_1=\cl{U}_0\cap \Delta (\cl{U})\subset \cl{U}_0.$\\
By the previous proposition we have: $tr(RQ_1R^*)=0\Rightarrow tr(R^*RQ_1)=0\Rightarrow RQ_1=0\Rightarrow R=R(I-Q_1).$\\
We conclude that $R_1(\Delta (\cl{U}))\subset \Delta (\mathcal{U})(I-Q_1).$

\medskip

For part (ii), observe that if $K\in \Delta (\cl{U})\cap \cl{K}$
then $K$ can be approximated in the norm topology by sums of rank
1 operators in $\Delta (\cl{U}):$ Proposition 3.4 in \cite{kt}. $ \qquad
\Box$

\begin{remark} We will see below that if  $\cl{U}$ is a strongly reflexive masa bimodule then $[R_1(\Delta (\cl{U}))]^{-w^*}=\Delta (\cl{U})(I-Q_1).$ This is not true in general. For example take $\cl{U}$ to be a    $TRO$ which is not strongly reflexive. Then $[R_1(\Delta (\cl{U}))]^{-w^*}$ is strictly contained in $\Delta (\cl{U})(I-Q_1)=\cl{U}.$
\end{remark}

\begin{proposition}
$\Delta(\cl{U})\subset (R_1(\cl{U}_0)^*)^0.$

\end{proposition}
\textbf{Proof}

\medskip

Let $T\in R_1(\mathcal{U}_0).$ Then as in Proposition 6.1 we have $T=L+M$ where  $$L\in R_1(\Delta(\mathcal{U}))\; \text{and}\; M\in [R_1(\phi (P_n)\cl{U}P_n^\bot):n\in \mathbb{N}]^{-\|\;\|_1}\subset \cl{U}_0.$$

\medskip

So $L=T-M\in \cl{U}_0\cap R_1(\Delta(\mathcal{U})). $

\medskip

Using Proposition 6.3,\;$\cl{U}_0\cap R_1(\Delta(\mathcal{U}))\subset \cl{U}_0\cap \Delta (\cl{U})(I-Q_1)$ which vanishes by Theorem 3.3 so $L=0$ and hence $T=M.$

\medskip

We conclude that \begin{equation}\label{4.1}R_1(\mathcal{U}_0)\subset [R_1(\phi (P_n)\cl{U}P_n^\bot):n\in \mathbb{N}]^{-\|\;\|_1}.\end{equation}

\medskip

Let $A\in \Delta (\cl{U}).$ We want to show that $tr(A^*R)=0$ for every $R\in R_1(\mathcal{U}_0).$

\medskip

Using $(6.1)$ it suffices to show that $tr(A^*R)=0$ for every $R\in R_1(\phi (P_n)\cl{U}P_n^\bot),$ and $n\in \mathbb{N}.$

\medskip

Let $R$ a rank 1 operator such that $R=\phi (P_n)RP_n^\bot$ then
 \begin{align*}tr(A^*R)=&tr(A^*\phi (P_n)RP_n^\bot)=tr(P_n^\bot A^*\phi (P_n)R)\\=&tr((\phi (P_n)AP_n^\bot)^*R)=tr(0R)=0.  \qquad  \Box
\end{align*}

\bigskip

Let $P\in \cl{S}_{1,\phi }.$ We suppose that $\vee \{\phi (L): L\in \cl{S}_{1,\phi },\phi (L)<\phi (P) \}<\phi (P).$\\
Since $\cl{S}_{2,\phi }$ is join complete there exists $P_0\in \cl{S}_{1,\phi }$ such that $$\phi (P_0)=\vee \{\phi (L): L\in \cl{S}_{1,\phi },\phi (L)<\phi (P) \}.$$

\medskip

We call the projection $P-P_0$ \textbf{an\;atom} of $\cl{U}$ and we denote the projection $\phi (P)-\phi (P_0)$ by $\delta (P-P_0).$

\begin{proposition}Let $F$ be an atom of $\cl{U}.$\\
i)The projection $F$ is minimal in the algebra $(\cl{S}_{1,\phi })^{\prime\prime}.$\\
ii)The projection $\chi (I)\delta (F)$ is minimal in the algebra $\chi (I)(\cl{S}_{2,\phi })^{\prime\prime}.$\\
iii)$\chi (I)\delta  (F)B(H_1,H_2)F\subset \Delta (\cl{U}).$\\
iv)$\chi (I)^\bot  \delta  (F)B(H_1,H_2)F\subset \cl{U}_0.$
\end{proposition}
\textbf{Proof}

\medskip

i) Let $P,P_0\in \cl{S}_{1,\phi }$ be such that $\phi (P_0)=\vee \{\phi (L): L\in \cl{S}_{1,\phi },\phi (L)<\phi (P) \}<\phi (P)$ and $F=P-P_0.$\\
If $Q\in \cl{S}_{1,\phi }$ either $P\leq Q$ or $QP<P.$\\
If $P\leq Q$ then $QF=F.$\\
If $QP<P$ then (since $QP\in \cl{S}_{1,\phi }$ and $\phi$ is 1-1 on $\cl{S}_{1,\phi }$)\;$\phi(QP)<\phi (P)\Rightarrow \phi(QP)\leq \phi (P_0)\Rightarrow QP\leq P_0,$ so $QF=0.$

\medskip

We conclude that $QFB(H_1)F=FB(H_1)QF$ for all $Q\in \cl{S}_{1,\phi },$ therefore $FB(H_1)F\subset (\cl{S}_{1,\phi })^\prime,$ hence $F$ is a minimal projection in $(\cl{S}_{1,\phi })^{\prime\prime}.$

\medskip

ii)Since $P,P_0\in \cl{S}_{1,\phi }$ we have that $\phi (P)\Delta (\cl{U})=\Delta (\cl{U})P$ and  $\phi (P_0)\Delta (\cl{U})=\Delta (\cl{U})P_0$ hence $$\delta (F)\Delta (\cl{U})=\Delta (\cl{U})F\; \text{and\;so}\; \chi (I)\delta (F)=\chi (F).$$
Let $Q\in \cl{S}_{1,\phi }.$

\begin{equation}\label{1}If \qquad QF=0\; \text{then}\; \chi (I)\delta (F)\phi (Q)=0.\end{equation}
Indeed, $\delta (F)\Delta (\cl{U})=\Delta (\cl{U})F$ so $\delta (F)\Delta (\cl{U})Q=0$ so $\delta (F)\chi (Q)=0$ so $\chi (I)\delta (F)\phi (Q)=0.$

\begin{equation}\label{2}If \qquad QF=F\Rightarrow \chi (I)\delta (F)\phi (Q)=\chi (I)\delta (F).\end{equation}
Indeed, $\delta (F)\Delta (\cl{U})=\Delta (\cl{U})F$ \;so\; $\delta (F)\Delta (\cl{U})Q=\Delta (\cl{U})F$\; so\;  $\delta(F) \chi (Q)=\chi (F)$ \;so\; $\chi (I)\delta (F)\phi (Q)=\chi (I)\delta (F).$

\medskip

Using equations $(6.2),(6.3)$ as in (i) we have that $\chi (I)\delta (F)$ is a minimal projection in $\chi (I)(\cl{S}_{2,\phi })^{\prime\prime}.$

\medskip

iii)Let  $T\in B(H_1,H_2)$ and $Q\in \cl{S}_{1,\phi }.$\\
From equations $(6.2),(6.3)$ it follows that $\phi (Q)\chi (I)\delta (F)TF=\chi (I)\delta (F)TFQ,$\;
so \;$\chi (I)\delta (F)TF\in \Delta (\cl{U}).$

\medskip

iv)If $T\in \cl{U}$ then $\chi (I)^\bot T\in \cl{U}_0.$ Indeed, by Theorem 3.1 there exist $T_1\in \cl{U}_0,\; T_2\in \Delta (\cl{U})$ so that $T=T_1+T_2.$\\
But $T_2=\chi (I)T_2$ so $\chi (I)^\bot T=\chi (I)^\bot T_1\in \cl{U}_0.$

\medskip

Now it suffices to show that $\delta (F)B(H_1,H_2)F\subset \cl{U}.$\\
Let $T\in B(H_1,H_2)$ and $Q\in \cl{S}_{1,\phi }.$\\
If $FQ=0$ then $\phi (Q)^\bot \delta (F)TFQ=0.$\\
If $FQ=F$ then $P-P_0\leq Q$ hence $\delta (F)=\phi (P)-\phi (P_0)\leq \phi (P-P_0)\leq \phi (Q)$ so
$\phi (Q)^\bot \delta (F)TFQ=0.$\\
We conclude that $\delta  (F)TF\in \cl{U}. \qquad  \Box$

\begin{remark}\emph{There exists a simple example of a reflexive masa bimodule $\cl{U}$ so that $\delta (F)B(H_1,H_2)F\subset \cl{U}_0$ for any atom $F$ in $\cl{U}.$ \\(Take $\cl{U}$ to be the set of $3\times 3$ matrixes with zero diagonal.)\\
This is an example of the different behaviour of algebras and bimodules:\\
it is known that if $\cl{U}$ is a CSL algebra in a Hilbert space $H$ and $F$ is an atom in $\cl{U}$ then $FB(H)F\subset \Delta (\cl{U}).$}
\end{remark}
We thank Dr. I.Todorov for suggesting the \lq atomic decomposition\rq \;in the theorem below.
\begin{theorem}Let $\{F_n: n\in \mathbb{N}\}=\{F: F\; atom\; of\; \cl{U}\}.$ Then $$[R_1(\Delta (\cl{U}))]^{-w^*}=\sum_{n=1}^\infty \oplus \chi (I)\delta (F_n)B(H_1,H_2)F_n.$$
\end{theorem}
\textbf{Proof}

\medskip

By the previous proposition it follows that $$[R_1(\Delta (\cl{U}))]^{-w^*}\supset \sum_{n=1}^\infty \oplus \chi (I)\delta (F_n)B(H_1,H_2)F_n.$$

Let $R=x\otimes y^*\in \Delta (\cl{U}).$\\
For every $Q\in \cl{S}_{1,\phi }$ we have that $x\otimes (Qy)^*=(\phi (Q)x)\otimes y^*$ so \\$\phi (Q)x\neq0\Leftrightarrow Qy\neq0\Leftrightarrow \phi (Q)x=x\Leftrightarrow Qy=y.$

\medskip

Let $P=\wedge \{Q\in \cl{S}_{1,\phi} : Qy=y\},$ then $P\in \cl{S}_{1,\phi }.$\\
If $Q\in \cl{S}_{1,\phi }$ so that $\phi(Q)<\phi (P)$ then $\phi (Q)x=0.$\\
(If $\phi(Q)x=x$ then $Qy=y$ so $Q\geq P$).

\medskip

Let $P_0\in \cl{S}_{1,\phi }$ with $\phi (P_0)=\vee \{\phi (L): L\in \cl{S}_{1,\phi },\phi (L)<\phi (P) \}.$\\
We observe that $\phi (P_0)x=0$ and $\phi (P)x=x,$ hence $\phi (P_0)<\phi (P).$\\
We conclude that $F=P-P_0$ is an atom of $\cl{U}.$\\
The equalities $(P-P_0)y=y$ and $(\phi (P)-\phi (P_0))x=x$ imply
that\\ $R=\delta (F)RF.$ But $R=\chi (I)R$ so  $R=\chi (I)\delta (F)RF.$\\
The proof is complete. $\qquad \Box$

\bigskip

Every strongly reflexive TRO is a masa bimodule \cite{kt}.\\
So using the previous theorem we have a new proof of the following result in \cite{kt}.

\begin{corollary}If $\cl{M}$ is a strongly reflexive TRO, $\zeta =\mathrm{Map} (\cl{M})$ and \\$\{A_n: n\in \mathbb{N}\}=\{A: A \;atom \;of\; \cl{M}\},$ then $$\cl{M}=\sum_{n=1}^\infty \oplus \zeta (A_n)B(H_1,H_2)A_n.$$

\end{corollary}

Let $(P_n)\subset \cl{S}_{1,\phi }$ be a sequence such that $$\cl{U}=\{T\in\;B(H_1,H_2):\phi(P_n)^{\bot}TP_n=0\; \text{for all}\; n\in\;\mathbb{N} \}.$$ Let $V_n,U_n:B(H_1,H_2)\longrightarrow B(H_1,H_2),n\in \mathbb{N}$ be as in theorem 3.1.\\
By Theorem 6.8 $$[R_1(\Delta (\cl{U}))]^{-w^*}=\sum_{n=1}^\infty \oplus E_nB(H_1,H_2)F_n,$$ where $F_n$ atom of $\cl{U}$ and  $E_n=\chi (I)\delta (F_n)$ for all $n\in \mathbb{N}.$

Thus $[R_1(\Delta (\cl{U}))]^{-w^*}$ is the range of the contractive projection $D$ defined by  $$D:B(H_1,H_2)\longrightarrow B(H_1,H_2):D(T)=\sum_{n=1}^{\infty }E_nTF_n.$$

\begin{proposition}\label{45}Let $K\in \cl{K},$ then the sequence $(U_n(K))$ converges to $D(K)$ in norm.
\end{proposition}
\textbf{Proof}

\medskip

We observe that $(V_n|_{C_2})$ is a commuting sequence of
orthogonal projections in the Hilbert space $C_2.$

\medskip

Hence $(U_n|_{C_2})$ is a decreasing sequence of orthogonal
projections. Therefore if  $T\in C_2$ the sequence $(U_n(T))$
converges in the Hilbert-Schmidt norm $\|\cdot\|_2.$

\medskip

Let $K\in \cl{K}.$ Then for $\varepsilon >0$ there exist $K_{\varepsilon }\in C_2 $ such that   $\|K-K_{\varepsilon }\|<\frac{\varepsilon }{3}$ and $n_0\in \mathbb{N}$ such that $\|U_n(K_{\varepsilon })-U_m(K_{\varepsilon })\|_2<\frac{\varepsilon }{3}$ for every $n,m \geq n_0.$

\medskip

Then
\begin{align*}&\|U_n(K)-U_m(K)\| \\&\leq \|U_n(K)-U_n(K_{\varepsilon} )\|+\|U_n(K_{\varepsilon })-U_m(K_{\varepsilon })\|+\|U_m(K_{\varepsilon })-U_m(K)\| \\&\leq \|K-K_{\varepsilon }\|+\|U_n(K_{\varepsilon })-U_m(K_{\varepsilon })\|_2+\|K-K_{\varepsilon }\|<\varepsilon\end{align*}  for every $n,m\geq n_0.$

\medskip

Thus $(U_n(K))$ converges in norm. Let $D_1(K)=\|\cdot\|-\lim U_n(K).$

\medskip

Since $\phi (P_i)^\bot U_n(K)P_i=\phi (P_i)U_n(K)P_i^\bot =0$ for every $i=1,2,...n,$ the limit $D_1(K)$  belongs to the diagonal $\Delta (\cl{U}).$

\medskip

Since $\|U_n(K)\|\leq \|K\|$ for all  $n\in \mathbb{N}$,  $D_1$ is a contraction.

\medskip

We observe that if $K\in \Delta (\cl{U})\cap \cl{K}$ then $U_n(K)=K$ for all $n\in \mathbb{N}$ hence $D_1$ projects onto $\Delta (\cl{U})\cap \cl{K}.$

\medskip

Now $D_1|_{C_2}$ is the orthogonal projection onto  $\Delta
(\cl{U})\cap C_2,$ being the infimum of the sequence
$(U_n|_{C_2}).$

\medskip

We can also observe that $D|_{C_2}$ is an orthogonal projection in the Hilbert space $C_2.$

\medskip

If $T\in \Delta (\cl{U})\cap C_2$ then by Proposition 6.3\; $T=\sum_{n=1}^{\infty }E_nTF_n=D(T).$

\medskip

We conclude that $D|_{C_2}$ and $D_1|_{C_2}$ are both orthogonal projections onto $\Delta (\cl{U})\cap C_2,$ hence $D|_{C_2}=D_1|_{C_2}.$

\medskip

Since $C_2$ is norm dense in $\cl{K}$ and $D|_{\cl{K}},D_1$ are
norm continuous, $D|_{\cl{K}}=D_1.\qquad \Box$

\begin{proposition}Suppose that $\vee _nF_n=F.$ Then the sequence $(U_n(T)F)$ converges strongly to the operator $D(T)$ for every $T\in B(H_1,H_2).$
\end{proposition}
\textbf{Proof}

\medskip

First we observe that if $x\in F_m(H_1), m\in \mathbb{N}$, then
the operator $x\otimes x^*$ is in $(\cl{S}_{1,\phi })^\prime.$

\medskip

Indeed, let $y\in E_m(H_2),$ then $R=y\otimes x^*\in \Delta (\cl{U}).$\\It follows that $R^*R=\|y\|^2x\otimes x^*\in \Delta (\cl{U})^*\Delta (\cl{U})\subset (\cl{S}_{1,\phi })^\prime.$

\medskip

Let $T\in B(H_1,H_2)$ and $x\in F_m(H_1),m\in \mathbb{N}, \|x\|=1.$

\medskip

By Proposition 6.10 $$U_i(T x\otimes x^*)\stackrel{\|\cdot\|}{\rightarrow}D(T x\otimes x^*), i\rightarrow \infty ,$$
hence \begin{equation}\label{4.2}U_i(T x\otimes x^*)(x)\stackrel{\|\cdot\|}{\rightarrow}D(T x\otimes x^*)(x), i\rightarrow \infty \end{equation}
\begin{equation}\label{4.3}D(T x\otimes x^*)(x)=\sum_{n=1}^\infty E_n (T x\otimes x^*) F_n (x)=E_mT(x)\end{equation}
\begin{equation}\label{4.4}D(T)(x)=\sum_{n=1}^\infty E_n T F_n(x)=E_mT(x)\end{equation}
We have that $$V_i(T x\otimes  x^*)=\phi (P_i)\;(T\; x\otimes x^*)\; P_i+\phi (P_i)^\bot (T\; x\otimes x^*)\; P_i^\bot ,\;i\in \mathbb{N}$$
 since $x\otimes x^* \in (\cl{S}_{1,\phi })^\prime,$  $$V_i(T x\otimes  x^*)=(\phi (P_i)\;T\; P_i)\; (x\otimes x^*) +(\phi (P_i)^\bot T\; P_i^\bot)\; (x\otimes x^*) ,\;i\in \mathbb{N},$$
hence \begin{equation}\label{4.5}U_i(T x\otimes x^*)=U_i(T)x\otimes x^*\Rightarrow U_i(T x\otimes x^*)(x)=U_i(T)(x), i\in \mathbb{N}\end{equation}

\medskip

Using $(6.4),(6.5),(6.6),(6.7)$ $$U_i(T)(x)\stackrel{\|\cdot\|}{\rightarrow }D(T)(x),\;i\rightarrow \infty \;\text{for all}\; x\in [\bigcup _{m=1}^\infty F_m(H_1)] .$$

\medskip

Since the $U_i$ are contractions $U_i(T)(x)\stackrel{\|\cdot\|}{\rightarrow} D(T)(x),$ for all   $x\in  F(H_1)$

\medskip

Observe that $D(T)F=D(T).\qquad \Box$

\begin{remark}
The sequence $(U_n(T))$ has similar properties to the net of
finite diagonal sums in the case of nest algebras.(Propositions
6.10,6.11 are analogous to Propositions 4.3,4.4 in \cite{dav}.)
\end{remark}

\begin{theorem}Let $K\in \cl{U}$ be compact. Then there exist unique compact operators $K_1\in \cl{U}_0, K_2\in \Delta (\cl{U})$ such that $K=K_1+K_2.$ Moreover $K_2=D(K).$
\end{theorem}
\textbf{Proof}

\medskip

Let $K_2=D(K)$ and $K_1=K-K_2$ then\;$K_1=\lim (K-U_n(K))$ (Proposition 6.10).

\medskip

As in Theorem 3.1 $K-U_n(K)\in \cl{U}_0$ for all $n\in \mathbb{N}.$ Hence $K_1\in \cl{U}_0.$

\medskip

The decomposition $K=K_1+K_2$ in $\cl{U}_0+\Delta (\cl{U})\cap
\cl{K}$ is unique because by Proposition 6.3,  $\Delta (\cl{U})\cap
\cl{K}\subset \Delta (\cl{U})(I-Q_1),$ while by Theorem 3.3,
$\cl{U}=\cl{U}_0\oplus \Delta (\cl{U})(I-Q_1).\qquad \Box$

\begin{corollary}Let $F\in \cl{U}$ be a finite rank operator. Then there exist unique finite rank operators $F_1\in \cl{U}_0,F_2\in \Delta (\cl{U})$ such that $F=F_1+F_2.$ Moreover $rank F_2\leq rank F$ and $F_2=D(F).$
\end{corollary}
\noindent\textbf{Proof }

\medskip

It can be shown that for each $n\in \mathbb{N}$ we have
$\mathrm{rank} (U_n(F))\le\mathrm{rank} (F).$

\medskip

 Therefore if $F_2=\|\cdot\|$-$\lim U_n(F)$ then
$\mathrm{rank} (F_2)\le\mathrm{rank} (F)$ and $F_2=D(F).$

\medskip

Setting $F_1=F-F_2$ we obtain the desired decomposition. $\qquad \Box$

\begin{corollary} Let $K\in \cl{U}\cap C_p ,1\leq p < \infty .$ Then there exist unique operators $K_1\in \cl{U}_0\cap C_p ,K_2\in \Delta (\cl{U})\cap C_p$ such that $K=K_1+K_2.$ Moreover $\|K_2\|_p\leq \|K\|_p.$
\end{corollary}
\textbf{Proof}

\medskip

As in Theorem 6.13 $K=K_1+D(K)$ where $K_1\in \cl{U}_0.$

\medskip

We observe that $D(K)\in C_p$ and $\|D(K)\|_p\leq \|K\|_p.\qquad \Box$

\section{Decomposition of a strongly reflexive masa bimodule}
Let $\cl{U},\;\cl{U}_0,\;\Delta (\cl{U}),\;\phi ,\;\cl{D}_1,\;\cl{D}_2$ be as in section 3 and $\chi =\mathrm{Map} (\Delta (\cl{U}))$.

\medskip

We now assume that $\mathcal{U}$ is a strongly reflexive masa bimodule.

\begin{proposition} The space $\mathcal{U}_0$ is strongly reflexive.
\end{proposition}
\textbf{Proof}

\medskip

Let $T\in \cl{U},\;P\in \cl{S}_{1,\phi }.$ Since $\cl{U}$ is a strongly reflexive masa bimodule there exists a net $(R_i)\subset [R_1(\cl{U})]$ such that $R_i\stackrel{wot}\rightarrow T:$ Corollary 2.5 in \cite{eks}.\\
So we have that $\phi (P)R_iP^\bot \stackrel{wot}\rightarrow \phi (P)TP^\bot.$ Since $(\phi (P)R_iP^\bot )\subset [R_1(\cl{U}_0)]$ we conclude that $\phi (P)TP^\bot\in [R_1(\cl{U}_0)]^{-wot}.$\\
We proved that $\phi (P)\cl{U}P^\bot \subset [R_1(\cl{U}_0)]^{-wot}$ for all $P\in \cl{S}_{1,\phi }.$ Hence $\mathcal{U}_0=[R_1(\cl{U}_0)]^{-wot}.  \qquad \Box$

\begin{remark}
The diagonal of a strongly reflexive masa bimodule is not
nececcarily strongly reflexive. For example
if $\cl{U}$ is a nonatomic nest algebra, then $\Delta (\cl{U})$
does not contain rank 1 operators.\end{remark}

\begin{proposition}
$i)\;\mathcal{U}_0=\mathcal{U}\cap (R_1(\Delta(\mathcal{U}))^*)^0.$\\
$ii)\;\mathcal{U}_0\cap\Delta(\mathcal{U})=\Delta(\mathcal{U})\cap(R_1(\Delta(\mathcal{U}))^*)^0.$
\end{proposition}
\textbf{Proof}

\medskip

By Proposition 6.2 we have $\cl{U}_0 \subset (R_1(\Delta
(\cl{U}))^*)^0.$

It suffices to show that $\mathcal{U}\cap (R_1(\Delta(\mathcal{U}))^*)^0\subset \cl{U}_0.$

\medskip

Since $\cl{U}\cap (R_1(\Delta(\cl{U}))^*)^0$ is masa bimodule, as in Theorem 3.1  we can decompose it in the next sum: $$\cl{U}\cap (R_1(\Delta(\cl{U}))^*)^0=\cl{U}_0\cap (R_1(\Delta(\cl{U}))^*)^0+\Delta (\cl{U})\cap (R_1(\Delta(\cl{U}))^*)^0.$$
Now we must prove that $\Delta (\cl{U})\cap (R_1(\Delta(\cl{U}))^*)^0\subset \cl{U}_0.$

\medskip

Using Theorem 2.2, there exist projections $P_1\in \cl{D}_1, P_2\in \cl{D}_2$ such that $[R_1(\Delta(\cl{U}))]^{-w^*}=P_2\Delta (\cl{U})P_1$ and $\Delta (\cl{U})\cap (R_1(\Delta(\cl{U}))^*)^0=P_2^\bot \Delta (\cl{U})P_1^\bot .$

\medskip

Let $T\in \Delta (\cl{U})\cap (R_1(\Delta(\cl{U}))^*)^0.$\\
Since $\cl{U}$ is a strongly reflexive masa bimodule there exists a net $(R_i)\subset [R_1(\cl{U})]$ such that $R_i\stackrel{wot}\rightarrow T$ \cite{eks}.\\
By Proposition 6.1 there exist $M_i \in [R_1(\Delta (\cl{U}))],
L_i\in \cl{U}_0$ such that  $R_i=M_i+L_i.$

\medskip

 Thus  $M_i+L_i \stackrel {wot}{\rightarrow }T$ so $P_2^\bot M_iP_1^\bot +P_2^\bot L_iP_1^\bot \stackrel {wot}{\rightarrow }P_2^\bot TP_1^\bot $ and thus $P_2^\bot L_iP_1^\bot  \stackrel {wot}{\rightarrow }T.$ It follows that $T\in \cl{U}_0. \qquad \Box$

\begin{theorem}\label{54}$\mathcal{U}=\mathcal{U}_0\oplus[R_1(\Delta(\mathcal{U}))]^{-w^*}$.\end{theorem}
\textbf{Proof}

\medskip

By Theorem 2.2, $$\Delta (\cl{U})=\Delta (\cl{U})\cap (R_1(\Delta(\mathcal{U})^*)^0\oplus [R_1(\Delta (\cl{U}))]^{-w^*}$$ so by Proposition 7.3 $\Delta (\cl{U})=\cl{U}_0\cap \Delta (\cl{U})+ [R_1(\Delta (\cl{U}))]^{-w^*}.$

\medskip

Since $\cl{U}=\cl{U}_0+\Delta (\cl{U})$ we have that $\mathcal{U}=\mathcal{U}_0+[R_1(\Delta(\mathcal{U}))]^{-w^*}.$

\medskip

By Proposition 6.3 and Theorem 3.3 the previous sum is direct.$ \qquad \Box$

\bigskip

Propositions 3.6 and 3.7 have the following consequences:

\begin{corollary}
i)The following are equivalent:

a) $ R_1(\Delta(\mathcal{U}))=0.$

b) $ \Delta(\mathcal{U})^*\Delta(\mathcal{U})\subset\mathcal{L}_1\cap
\cl{A}_1.$

c) $ \Delta(\mathcal{U})\Delta(\mathcal{U})^*\subset\mathcal{L}_2\cap
\cl{A}_2.$

ii)The following are equivalent:

a)$  \Delta(\mathcal{U})\;\; \text{is strongly reflexive.}$

b)  $\Delta(\mathcal{U})\;(\mathcal{L}_1\cap \cl{A}_1)=0.$

c)  $(\mathcal{L}_2\cap \cl{A}_2)\;\Delta(\mathcal{U})=0.$

\end{corollary}

Theorems 6.8, 7.4 and Corollary 5.3 give the following form of the decomposition of $\cl{U}$ when it is a strongly reflexive C.S.L. algebra.

\begin{corollary}If $\cl{S}$ is a completely distributive CSL in a Hilbert space $H$ and $\{A_n: n\in \mathbb{N}\}=\{A: A\; \text{atom\;of\;} \cl{S}\}$ then:
$$Alg(\cl{S})=Rad(Alg(\cl{S}))^{-w^*}\oplus \sum_{n=1}^\infty \oplus A_nB(H)A_n.$$
\end{corollary}

\bigskip

Recall the notation $[R_1(\Delta (\cl{U}))]^{-w^*}=\sum_{n=1}^{\infty }\oplus \chi (I)\delta (F_n)B(H_1,H_2)F_n,$\\ where $\{F_n: n\in \mathbb{N}\}=\{F: F\; atom\; of\; \cl{U}\}$ and $$D:B(H_1,H_2)\longrightarrow B(H_1,H_2):D(T)=\sum_{n=1}^{\infty }\chi (I)\delta (F_n)TF_n.$$

\begin{proposition}Let $\theta  :\cl{U}\rightarrow \cl{U}$ be the projection onto $[R_1(\Delta(\mathcal{U}))]^{-w^*}$ defined by the decomposition in Theorem 7.4. Then $\theta  =D|_{\cl{U}}.$\end{proposition}

\textbf{Proof}

\medskip

Since $\cl{U}$ decomposes as the direct sum of the masa bimodules
$\cl{U}_0$ and $[R_1(\Delta(\mathcal{U}))]^{-w^*}$, the map
$\theta $ is a masa bimodule map: $$\theta  (D_2TD_1)=D_2\theta
(T)D_1$$ for every $T\in \cl{U}, D_1\in \cl{D}_1, D_2\in
\cl{D}_2.$

\medskip

Hence if $T\in \cl{U}:$ \begin{align*}\theta  (T)=\sum_{n=1}^\infty \chi (I)\delta (F_n)\theta
(T)F_n&=\sum_{n=1}^\infty \theta  (\chi (I)\delta (F_n)TF_n)\\&=\sum_{n=1}^\infty
\chi (I)\delta (F_n)TF_n=D(T). \qquad \Box\end{align*}

\begin{proposition}$\cl{U}_0=\{T\in \cl{U}: \chi (I)\delta (F)TF=0 \text{\;for\;every\;atom\;F\;of}\;\cl{U} \}.$
\end{proposition}
\textbf{Proof}

\medskip

Let $F$ be an atom of $\cl{U}.$\\
If $P\in \cl{S}_{1,\phi },$ as in Proposition 6.6 either $PF=F\Rightarrow P^\bot F=0$ or $PF=0\Rightarrow \chi (I)\delta (F)\phi (P)=0.$\\
So $\chi (I)\delta (F)\phi (P)TP^\bot F=0$ for all $P\in \cl{S}_{1,\phi }$ and $T\in \cl{U},$
thus $\chi (I)\delta (F)\cl{U}_0F=0$ for every atom $F.$\\
It follows that $\cl{U}_0\subset \{T\in \cl{U}: \chi (I)\delta (F)TF=0, \text{for\;every\;atom\;F\;in}\;\cl{U} \}.$

\medskip

For the converse, let $T\in \cl{U}: \chi (I)\delta (F)TF=0$ for every atom $F$ in $\cl{U}.$\\
By the previous proposition $D(T)=0,$ hence $T\in \cl{U}_0.  \qquad  \Box$

\bigskip

It is known that the linear span of the rank 1 operators in a strongly reflexive masa bimodule is wot dense in the module.

This is not true generally  for the ultraweak topology \cite{eks}.

For the previous problem we have the next equivalence in proposition 7.10.

Firstly, we need the following lemma.

\begin{lemma}If $\cl{U}$ is a reflexive masa bimodule (not necessarily strongly reflexive) then: $$[R_1(\cl{U})]^{-w^*}=[R_1(\cl{U}_0)]^{-w^*}\oplus [R_1(\Delta (\cl{U}))]^{-w^*}.$$
\end{lemma}
\textbf{Proof}

\medskip

Since $R_1(\Delta (\cl{U}))\subset \Delta (\cl{U})(I-Q_1)$ (Proposition 6.3), by Theorem 3.3 the previous sum is direct.

Clearly $$[R_1(\cl{U})]^{-w^*}\supset [R_1(\cl{U}_0)]^{-w^*}\oplus [R_1(\Delta (\cl{U}))]^{-w^*}.$$

For the converse, let $T\in [R_1(\cl{U})]^{-w^*}.$

There is a net $(R_i)\subset [R_1(\cl{U})]$ with $R_i \stackrel{w^*}{\longrightarrow}T.$

\medskip

As in Proposition 6.1, we may decompose $R_i=L_i+M_i$ where \\$L_i\in [R_1(\cl{U}_0)]^{-\|\cdot\|_1}$ and $M_i\in [R_1(\Delta (\cl{U}))]$ for all i.

\medskip

Since $M_i=D(R_i)$ (Corollary 6.14) and $D$ is $w^*-$continuous,\\ we have $M_i \stackrel{w^*}{\longrightarrow}M \in [R_1(\Delta (\cl{U}))]^{-w^*}.$

\medskip

So $L_i=R_i-M_i\stackrel{w^*}{\longrightarrow}T-M=L \in [R_1(\cl{U}_0)]^{-w^*}.$

\medskip

Thus $T=L+M \in [R_1(\cl{U}_0)]^{-w^*}\oplus [R_1(\Delta (\cl{U}))]^{-w^*}. \qquad \Box$

\begin{proposition}If $\cl{U}$ is a strongly reflexive masa bimodule, then: $$\cl{U}=[R_1(\cl{U})]^{-w^*}\Leftrightarrow \cl{U}_0=[R_1(\cl{U}_0)]^{-w^*}.$$
\end{proposition}
\textbf{Proof}

\medskip

Suppose $\cl{U}=[R_1(\cl{U})]^{-w^*}.$ Then by Theorem 7.4 we have \\$\mathcal{U}=\mathcal{U}_0\oplus[R_1(\Delta(\mathcal{U}))]^{-w^*}.$

\medskip

It follows from the previous lemma that $\cl{U}_0=[R_1(\cl{U}_0)]^{-w^*}.$

\medskip

If conversely $\cl{U}_0=[R_1(\cl{U}_0)]^{-w^*}$ then again by Theorem 7.4 $$\mathcal{U}=\mathcal{U}_0\oplus[R_1(\Delta(\mathcal{U}))]^{-w^*}=[R_1(\cl{U}_0)]^{-w^*}\oplus [R_1(\Delta (\cl{U}))]^{-w^*}=[R_1(\cl{U})]^{-w^*}$$ by Lemma 7.9.$\qquad  \Box$

\bigskip

{\em Acknowledgement:} I would like to express appreciation to
Prof. A. Katavolos for his helpful comments and suggestions during
the preparation of this work, which is part of my doctoral thesis.
I wish to thank him for the proof of Proposition 6.10.

I am also indebted to Dr. I.Todorov for helpful discussions and
important suggestions that led to a substantial improvement of
this paper.

This research was partly supported by Special Account Research
Grant No. 70/3/7463 of the University of Athens.

GEORGE ELEFTHERAKIS\\
 Department of Mathematics,\\ University of
Athens, \\ Panepistimioupolis 157 84\\ Athens Greece

e-mail: gelefth@math.uoa.gr

\end{document}